\documentclass{article}

\usepackage[T1]{fontenc}
\usepackage{amsmath}
\usepackage{amsfonts}
\usepackage{mathrsfs}
\usepackage{yfonts}
\usepackage[dvips]{graphicx}

\newtheorem{theorem}{Theorem}
\newtheorem{lemma}[theorem]{Lemma}
\newtheorem{definition}[theorem]{Definition}
\newtheorem{example}[theorem]{Example}

\newtheorem{conjecture}[theorem]{Conjecture}

\newenvironment{proof}{\noindent\textsc{Proof: }}{\hspace{\stretch{1}}$\square$\medskip}

\title{Transitive graphs in counterexamples to 
   Karp's conjecture}
\author{Alexander Engström\footnote{
Research supported by ETH and Swiss National Science Foundation Grant PP002-102738/1. E-mail: engstroa@inf.ethz.ch 
}
\\ \\
\emph{Institute of Theoretical Computer Science,} \\ 
\emph{ETH Z\"urich, CH-8092 Z\"urich, Switzerland}}

\begin{document}

\maketitle

\noindent\textbf{Abstract:} Karp conjectured
that all nontrivial monotone graph properties 
are evasive. This was proved for $n$ a prime 
power, and $n=6$, where $n$ is the number of 
graph vertices, by Kahn, Saks, and Sturtevant.
We give a complete description of which 
transitive graphs are contained in a possible 
counterexample when $n=10$. 

\section{Introduction}

The notion of evasiveness came from the study
of argumented complexity, but with Kahn, Saks
and Sturtevants influential paper \cite{KSS} 
it was incorporated in combinatorial and
algebraic topology. A graph property is a
partition of the unlabeled graphs into two
classes, those with and those without the 
property. If the property is preserved under the
removal of edges it is monotone. Let us fix
a monotone graph property and the number of
graph vertices, then we can create a simplicial 
complex with the graph edges as vertices and
the graphs with the property as simplices.
\begin{example}
Let the graph property be planarity and use 
5 graph vertices. All but the complete graph are
planar, so the faces of the simplicial complex
are all proper subsets of the graph edge set.
\end{example}
\begin{definition}
A simplicial complex $\Delta$ is nonevasive
if it is a point, or if there is a vertex $v$
of $\Delta$ such that both the link 
$\mathrm{lk}_\Delta(v)=\{\sigma\in \Delta
\mid v \not\in \Delta, \sigma\cup \{v\}\in\Delta
\}$ and the deletion $\mathrm{dl}_\Delta(v)=
\{\sigma\in \Delta \mid v \not\in \Delta\}$ 
are nonevasive.
\end{definition}
A simplicial complex is evasive if it is not
nonevasive and a monotone graph property on a 
certain number of vertices is evasive or 
nonevasive dependent on its simplicial complex.
A trivial graph properties include all or
none graphs.
Now we can state the famous conjecture by Karp.
\begin{conjecture}
All nontrivial monotone graph properties are
evasive.
\end{conjecture}
Kahn et al \cite{KSS} used fixed-point theorems
by Oliver \cite{O} and Smith \cite{S} and
the implications (see \cite{B,KSS})
\[\mathrm{nonevasive}\Rightarrow 
  \mathrm{collapsible}\Rightarrow
  \mathrm{contractible}\Rightarrow
  \mathbb{Z}-\mathrm{acyclic}\Rightarrow
  \mathbb{Z}_p-\mathrm{acyclic}
\]
to prove the conjecture when the number of
graph vertices is a prime power or 6. The
goal of this paper is a characterization of
which vertex transitive graphs are in a possible
counterexample of the conjecture for graph
properties on 10 vertices. Our method is the
topological and we will use these results
by Oliver~\cite{O}:
\begin{theorem}\label{theorem:1}
If $\Gamma' \lhd \Gamma$, $\Gamma / \Gamma'$
is cyclic, $\Gamma'$ is of $p$ prime
power order, and $\Delta$ is 
$\mathbb{Z}_p$--acyclic, then 
$\chi(\Delta^\Gamma)=1$.
\end{theorem}
\begin{theorem}\label{theorem:2}
If $\Gamma '' \lhd \Gamma' \lhd \Gamma$, 
$\Gamma' / \Gamma''$ is cyclic, 
$\Gamma''$ is of $p$ prime
power order,
$\Gamma / \Gamma'$ is of $q$ prime
power order,
 and $\Delta$ is 
$\mathbb{Z}_p$--acyclic, then 
$\chi(\Delta^\Gamma)\equiv 1\: \mathrm{mod}\: q$.
\end{theorem}
\begin{example}
To illustrate the method, we prove the
conjecture when there are five vertices.
Label the vertices $1,2,3,4,5$.
Assume that the conjecture is false and create
a simplicial complex $\Delta$ from that
graph property. The graph property is nontrivial,
so $\Delta$ is neither empty nor a full simplex.
Since $\Delta$ is nonevasive it is 
$\mathbb{Z}$-acyclic and we can use 
Theorem~\ref{theorem:1}.
The action of the cyclic group $\Gamma=\langle 
(1\:2\:3\:4\:5) \rangle$ on the graph 
vertices, induces an action on the graph edges,
which is the same as the vertices of $\Delta$.
The abstract simplicial complex $\Delta^\Gamma$
has the minimal nonempty $\Gamma$-invariant
faces of $\Delta$ as vertex set, and a set of
vertices of $\Delta^\Gamma$ is a face of
$\Delta^\Gamma$ if their union is a face of 
$\Delta$. The minimal nonempty $\Gamma$-invariant
graphs are:
\begin{center}
\begin{tabular}{cc}
\includegraphics*[width=0.5\textwidth]
{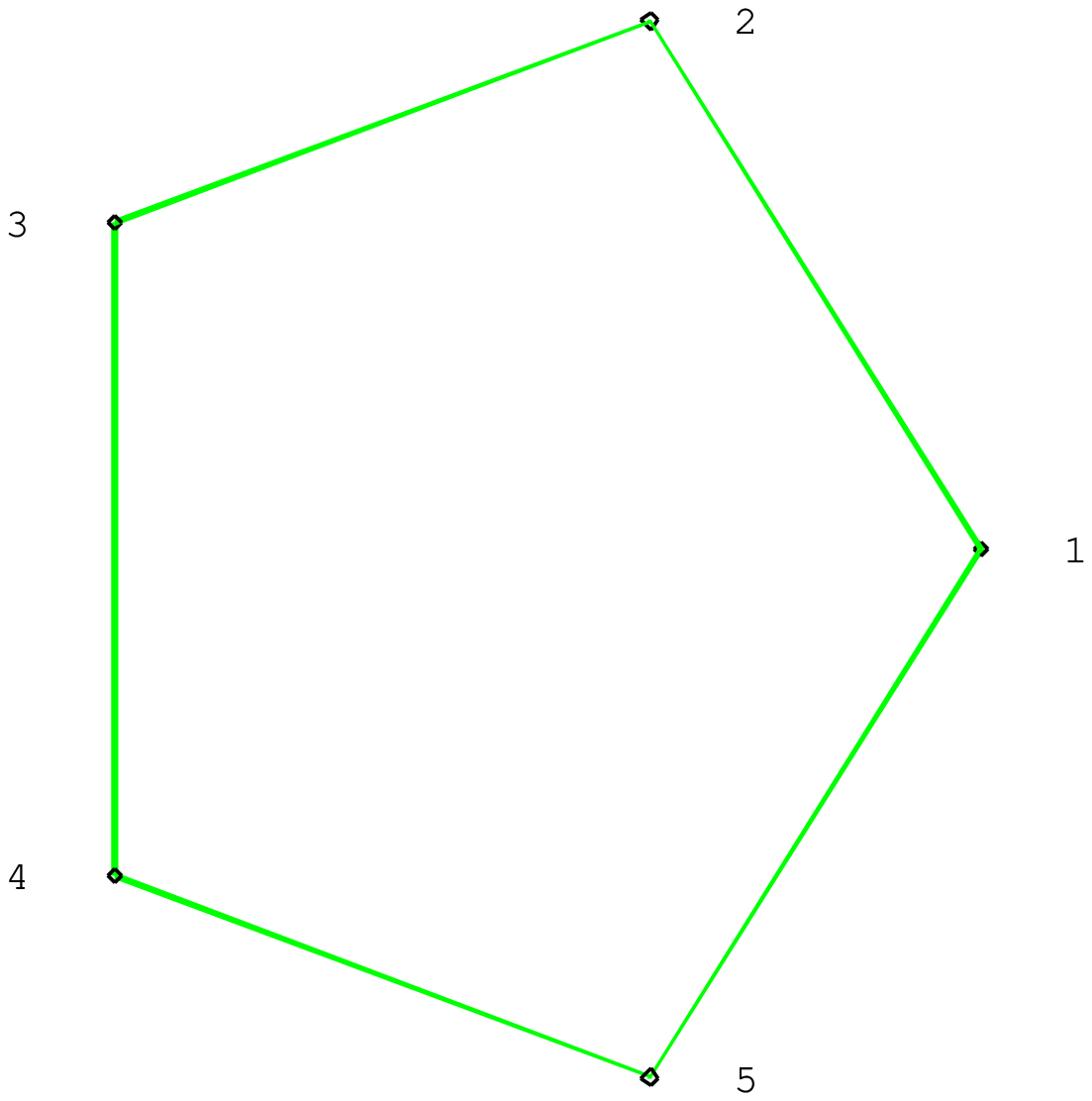} &
\includegraphics*[width=0.5\textwidth]
{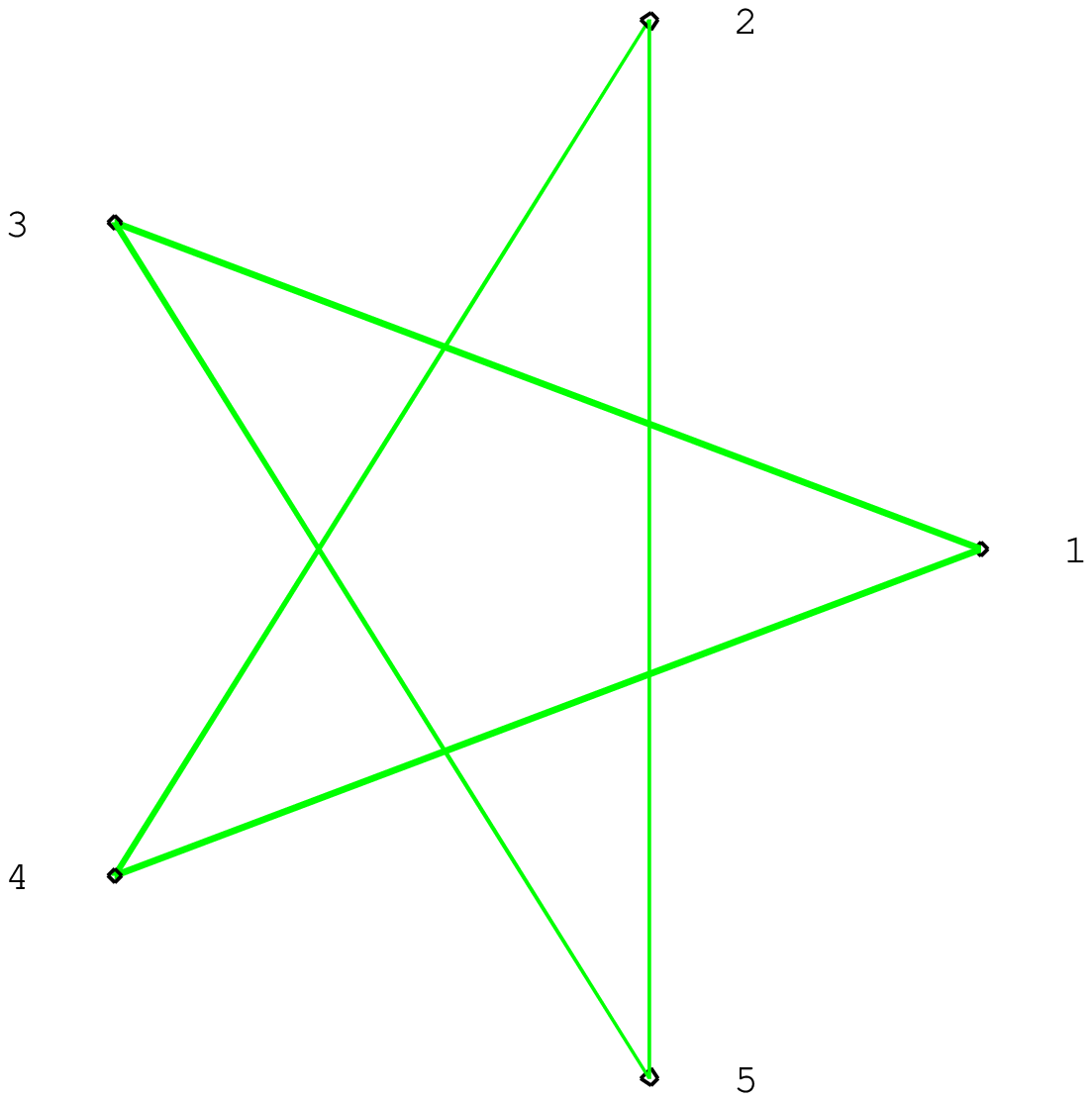}
\end{tabular}
\end{center}
If the five cycle is in $\Delta$ then 
$\Delta^\Gamma$ is two disjoint points, and
$\chi(\Delta^\Gamma)=2$. If not, then 
$\chi(\Delta^\Gamma)=0$. Using 
Theorem~\ref{theorem:1} with $\Gamma'$ as the
trivial group we get that 
$\chi(\Delta^\Gamma)=1$,
which is a contradiction. Hence the conjecture
is true for five vertices.
\end{example}

The plan of the rest of the paper is:
First we describe all transitive graphs on ten
vertices and their inclusion order. Then we
use Theorem \ref{theorem:1} and \ref{theorem:2}
to get conditions on which graphs can be in
a counterexample of the conjecture. The
inclusion order and conditions are investigated
and they give six different cases for which
transitive graphs are in a counterexample.
As an appendix we describe some computational
methods.

The reader who wants a more elaborate description
of the method herein could look at the original
work by Kahn et al \cite{KSS} and for example
Lutz \cite{L1,L2}. Nonevasiveness has poped up
outside the domain of graph properties, see
Kozlov \cite{Ko}, Kurzweil \cite{Ku}, and 
Welker \cite{W}.

\section{The transitive graphs}
\subsection{Cayley graphs}
There are 22 transitive graphs on 10 vertices,
and 20 of them are Cayley graphs \cite{McK,R}.
\subsubsection{Definitions}
\begin{definition}
Let $D$ be a subset of $\{1,2,3,4,5\}$. The 
vertex set of $G_D$ is $\{1,2,\ldots 10\}$, and
two vertices $i>j$ are adjacent if $(i-j)\in D$ 
or $(10+j-i)\in D$.
\end{definition}
Notice that $G_{ \{1,2,3,4,5\}}$ is the complete
graph $K_{10}$, and that the complement of $G_D$
is $G_{\{1,2,3,4,5\}\setminus D}$. Some of the
graphs are isomorphic.
\begin{center}
\begin{tabular}{cc}
\includegraphics*[width=0.5\textwidth]
{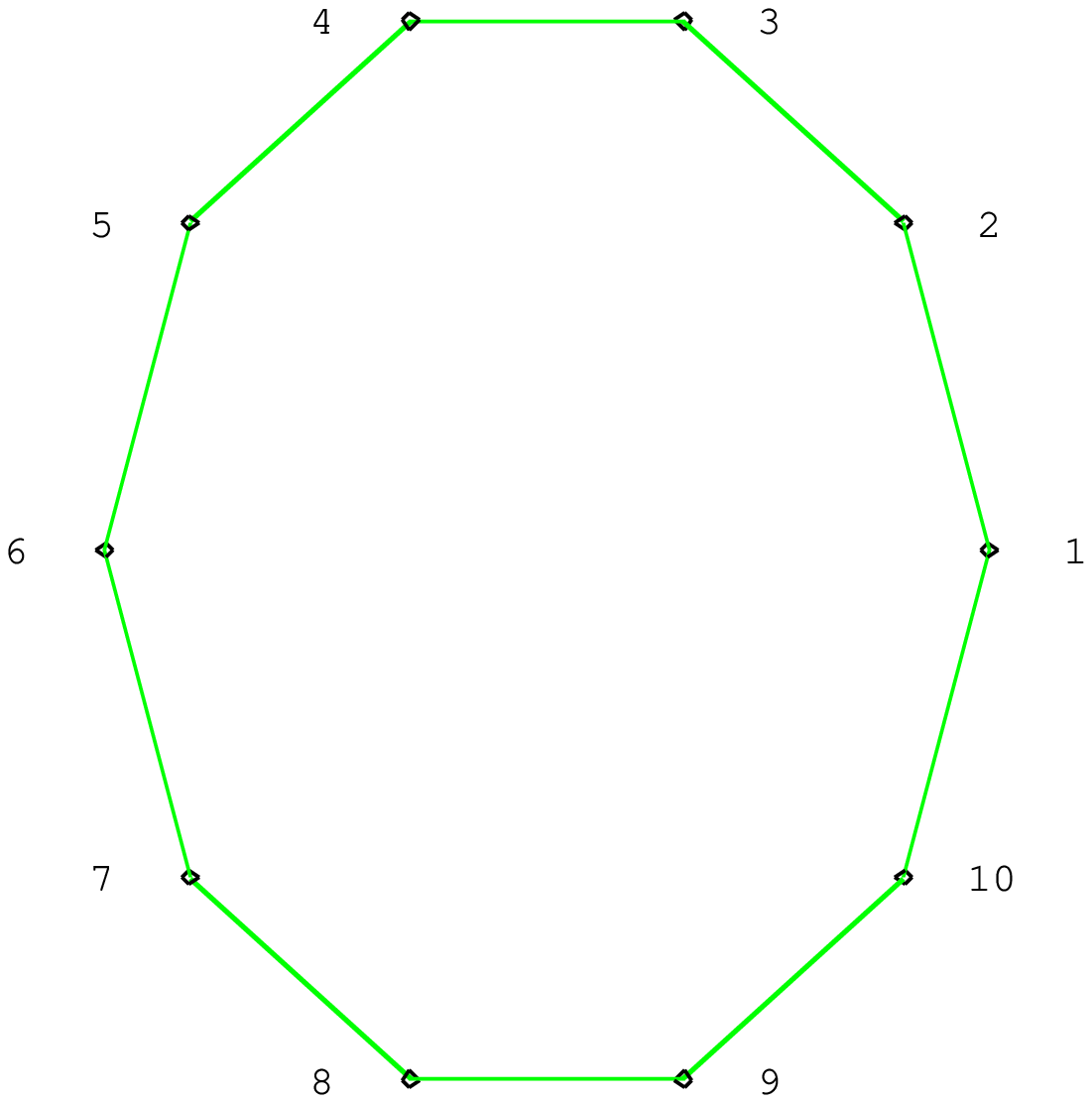} &
\includegraphics*[width=0.5\textwidth]
{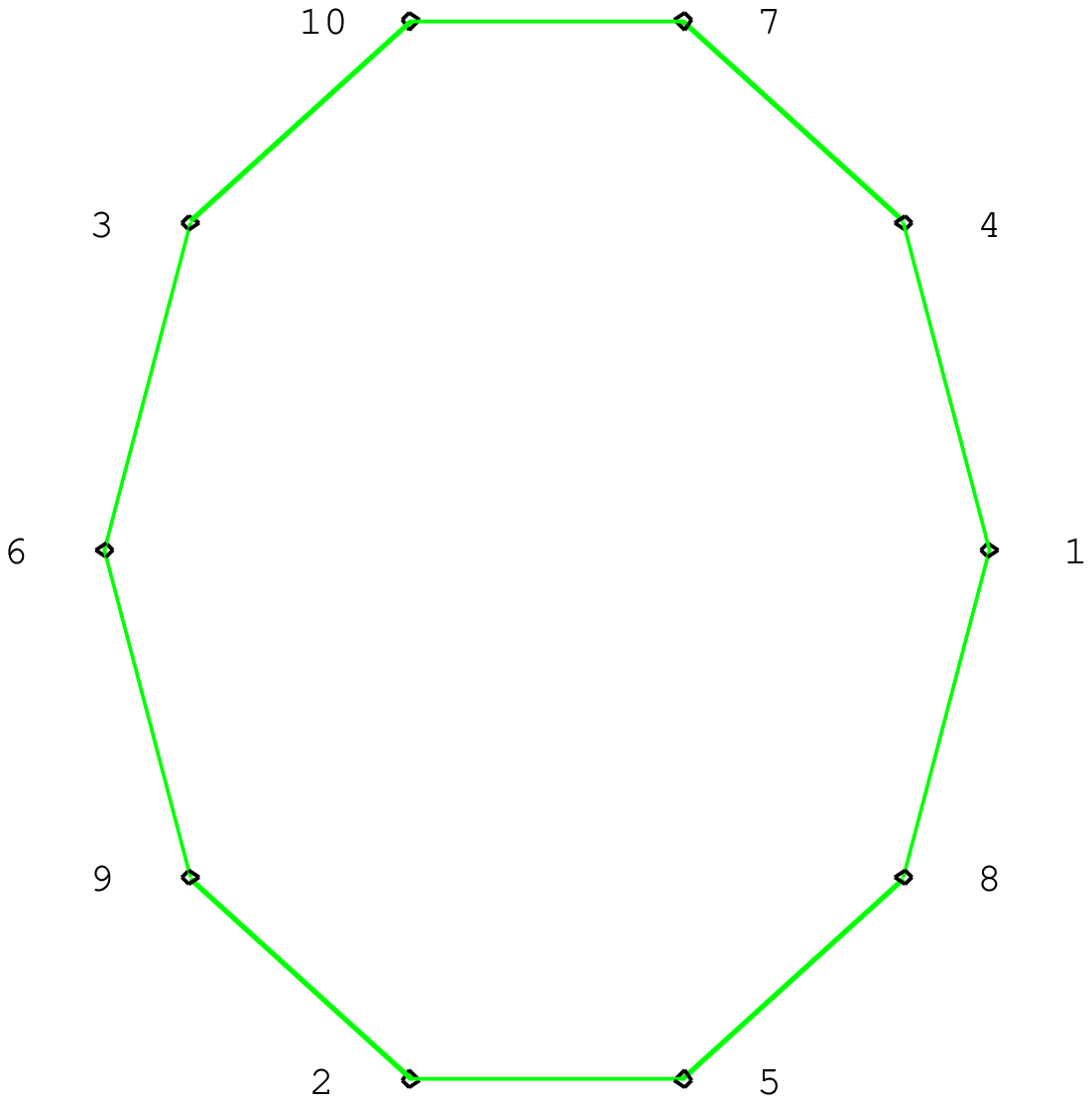}
\end{tabular}
$G_{\{1\}} \simeq G_{\{3\}}$ and 
$G_{\{2,3,4,5\}} \simeq G_{\{1,2,4,5\}}$.
\end{center}

\begin{center}
\begin{tabular}{cc}
\includegraphics*[width=0.5\textwidth]
{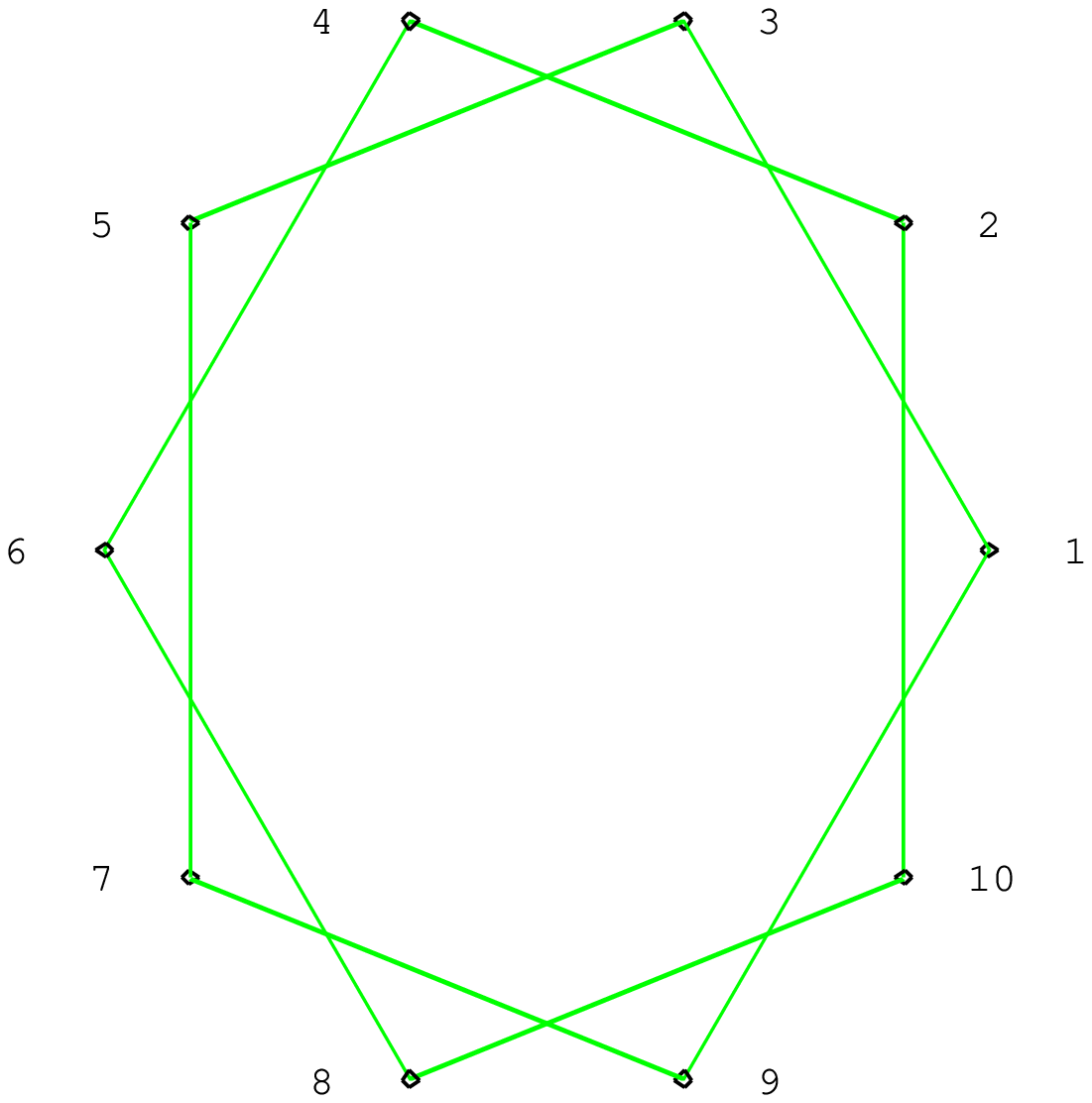} &
\includegraphics*[width=0.5\textwidth]
{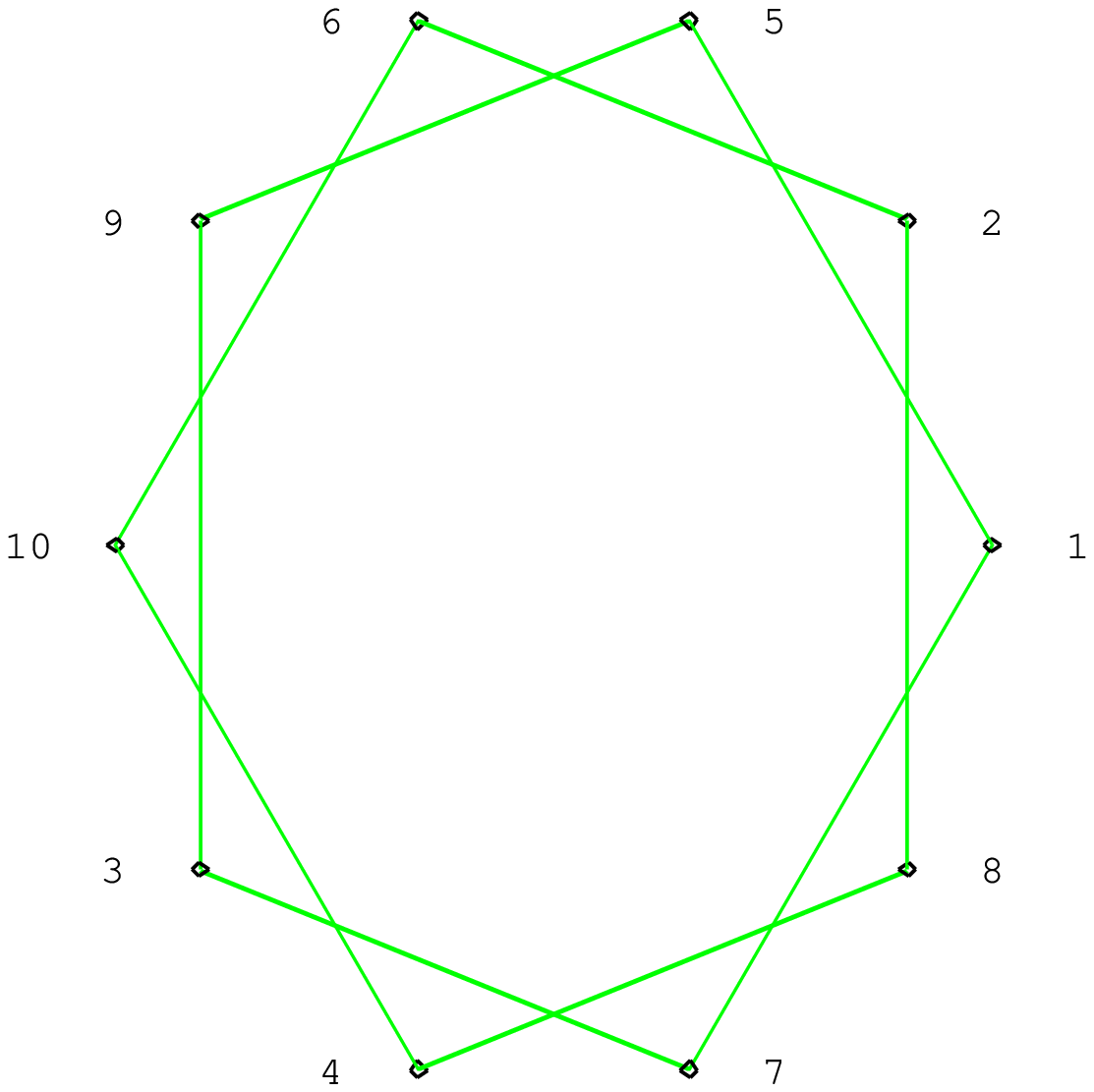}
\end{tabular}
$G_{\{2\}} \simeq G_{\{4\}}$ and 
$G_{\{1,2,3,5\}} \simeq G_{\{1,3,4,5\}}$.
\end{center}

\begin{center}
\begin{tabular}{cc}
\includegraphics*[width=0.5\textwidth]
{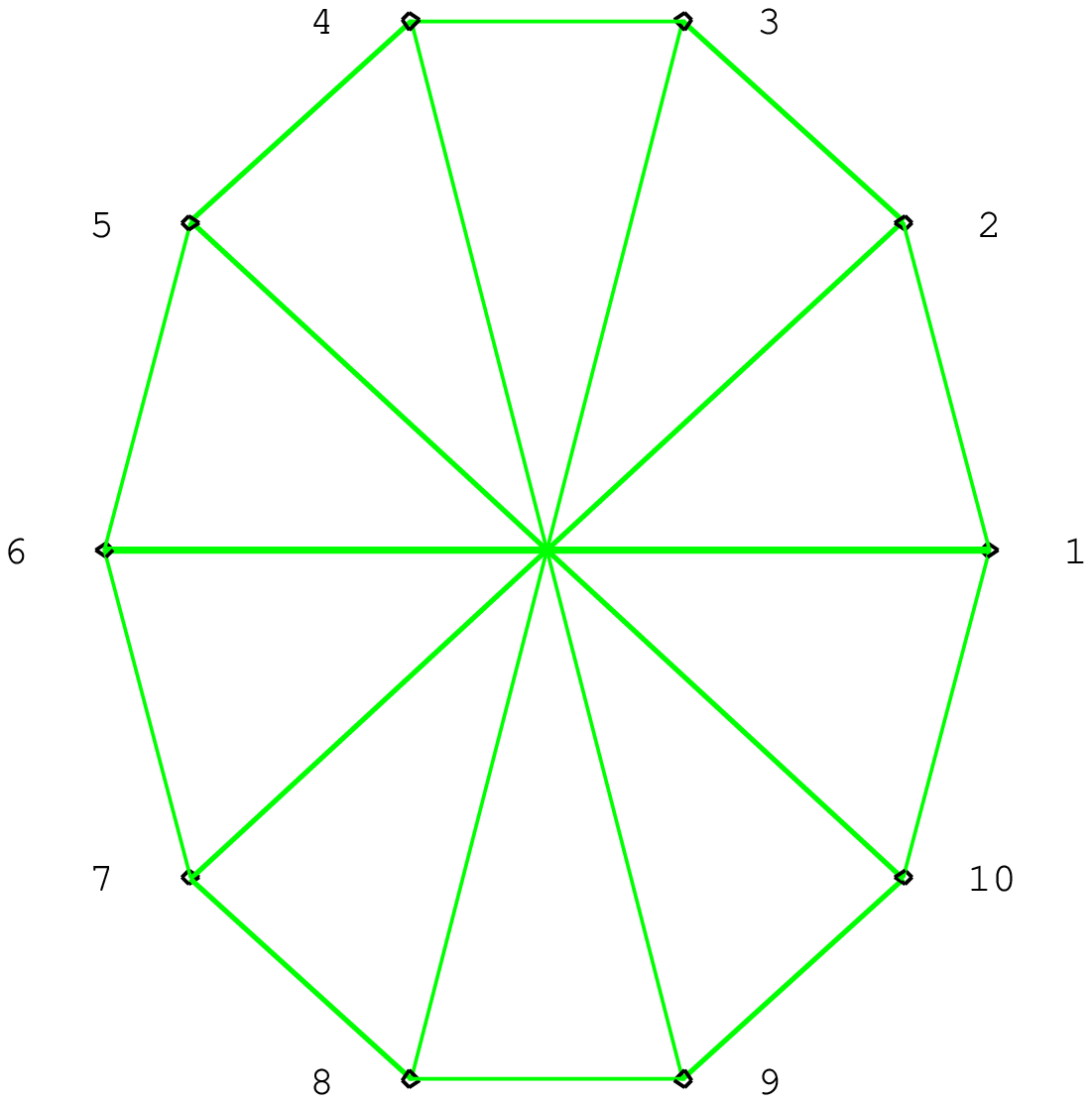} &
\includegraphics*[width=0.5\textwidth]
{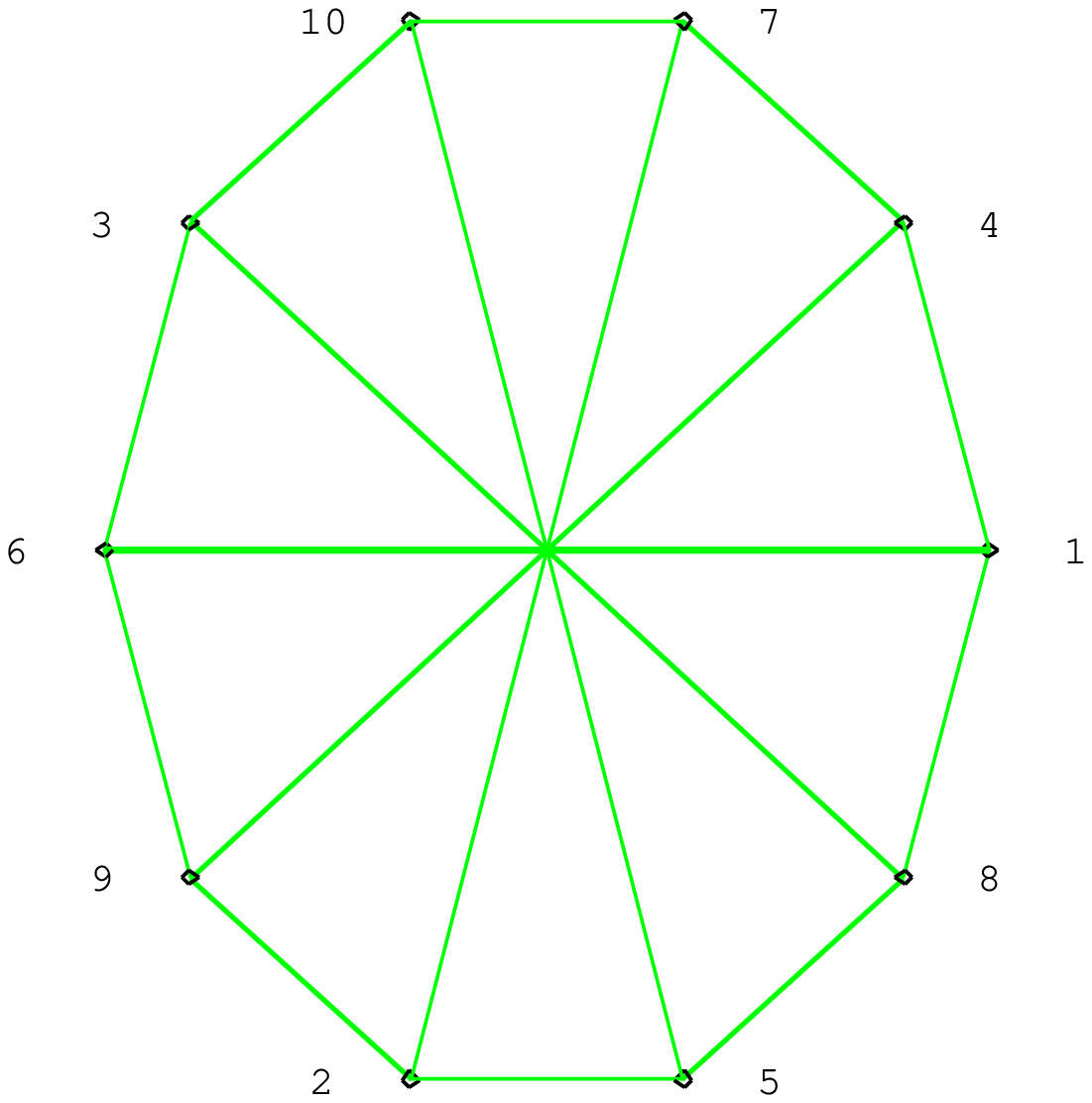}
\end{tabular}
$G_{\{1,5\}} \simeq G_{\{3,5\}}$ and 
$G_{\{1,2,4\}} \simeq G_{\{2,3,4\}}$.
\end{center}

\begin{center}
\begin{tabular}{cc}
\includegraphics*[width=0.5\textwidth]
{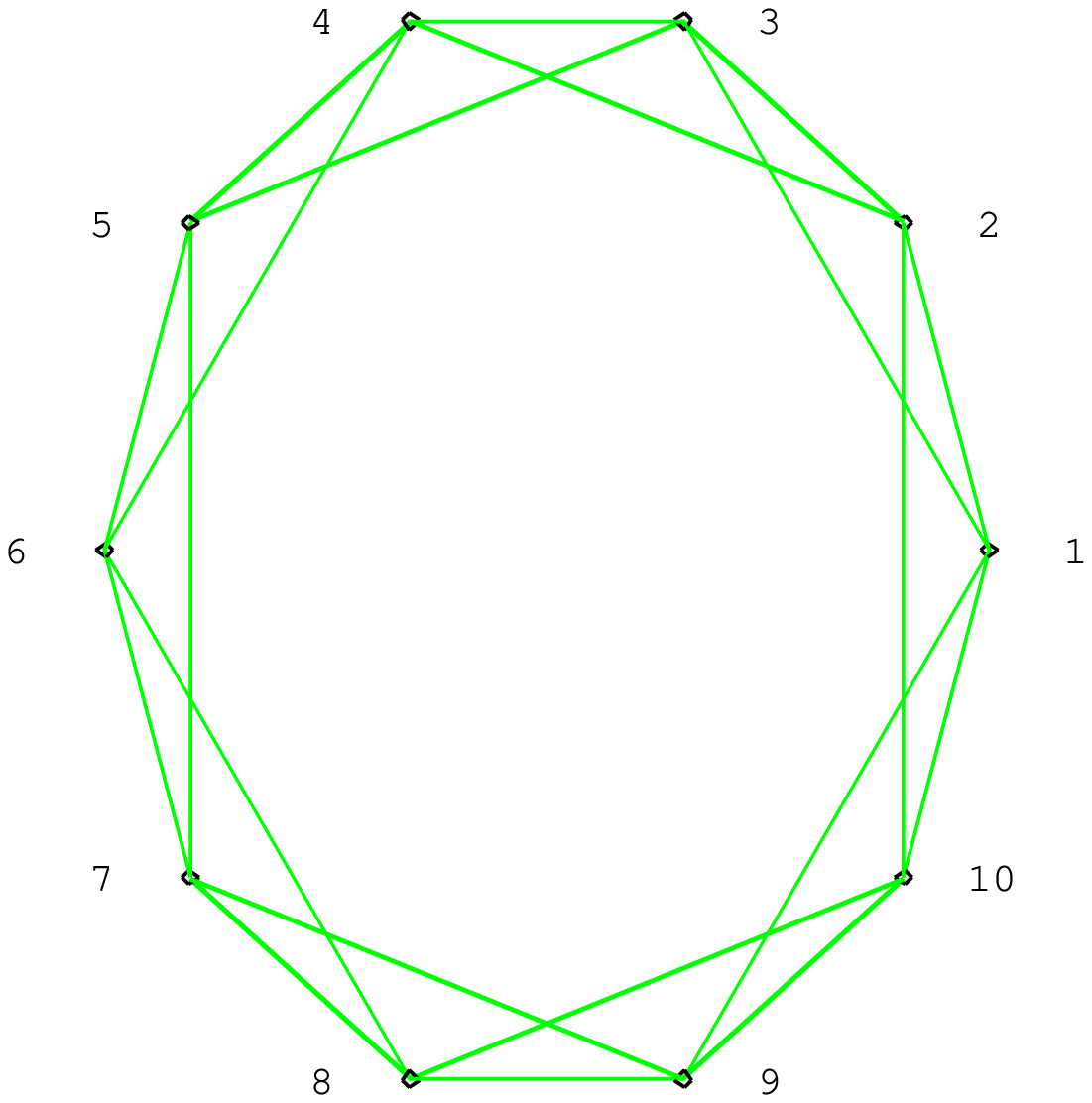} &
\includegraphics*[width=0.5\textwidth]
{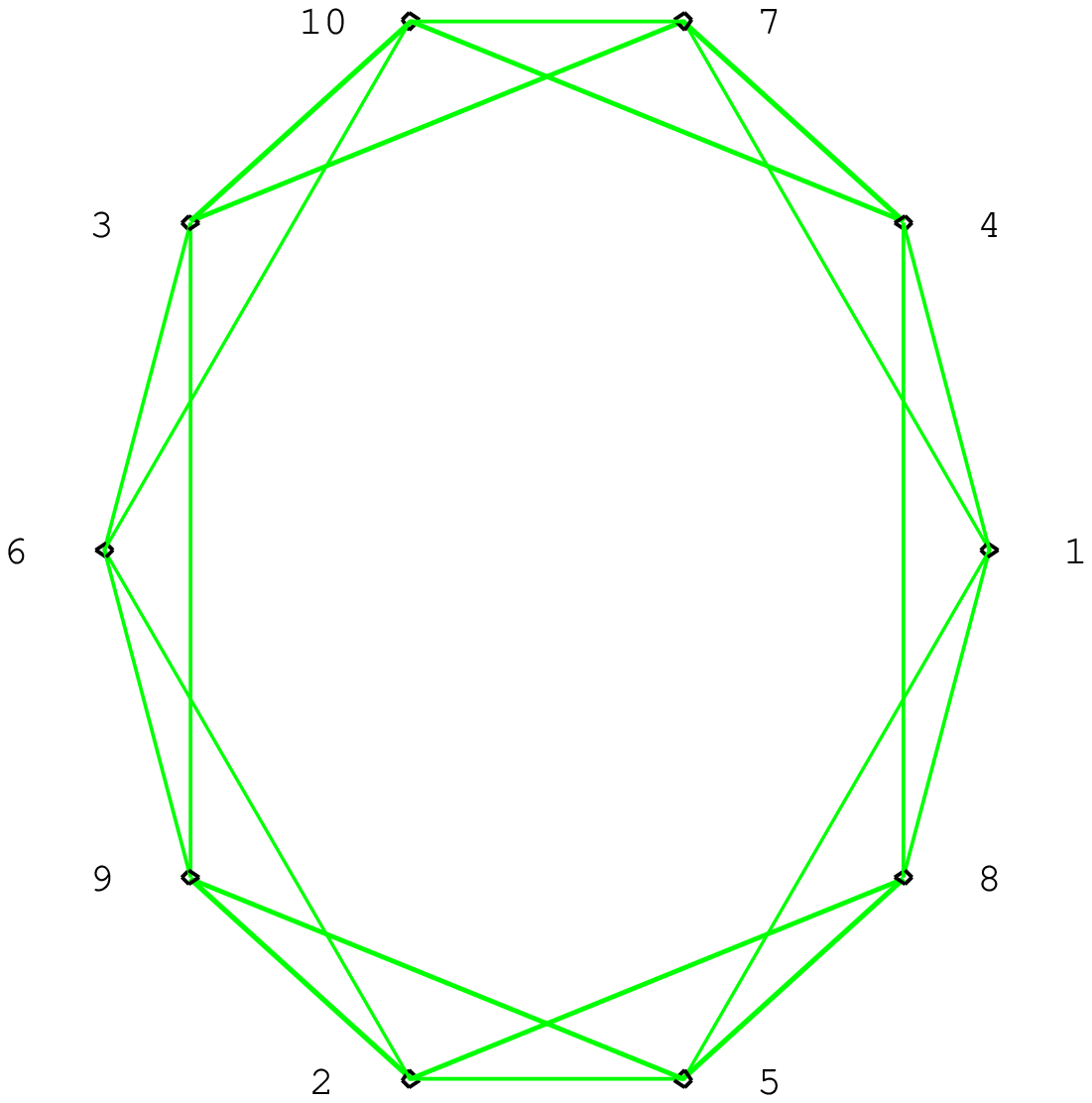}
\end{tabular}
$G_{\{1,2\}} \simeq G_{\{3,4\}}$ and 
$G_{\{1,2,5\}} \simeq G_{\{3,4,5\}}$.
\end{center}

\begin{center}
\begin{tabular}{cc}
\includegraphics*[width=0.5\textwidth]
{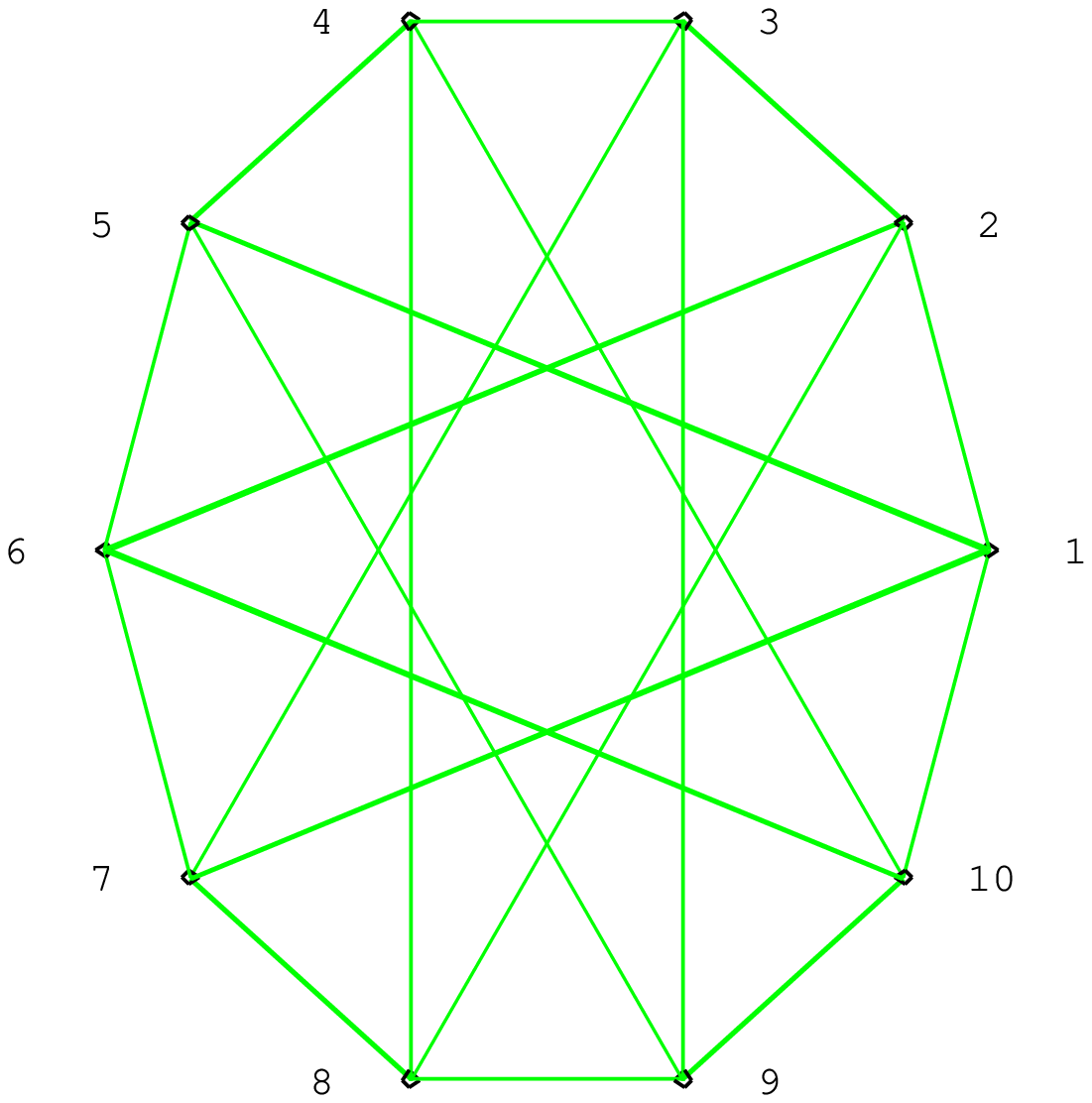} &
\includegraphics*[width=0.5\textwidth]
{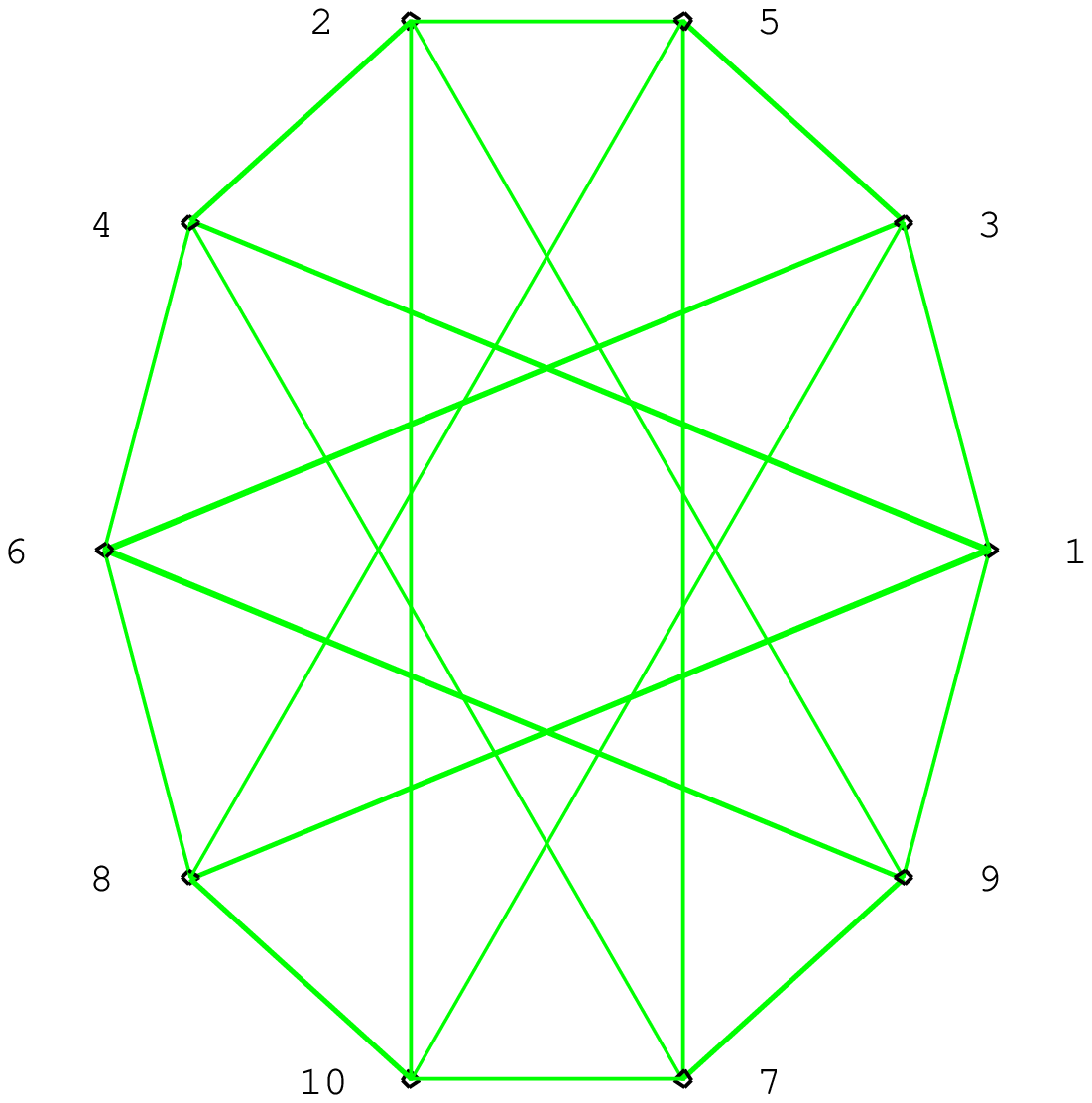}
\end{tabular}
$G_{\{1,4\}} \simeq G_{\{2,3\}}$ and 
$G_{\{1,4,5\}} \simeq G_{\{2,3,5\}}$.
\end{center}

\begin{center}
\begin{tabular}{cc}
\includegraphics*[width=0.5\textwidth]
{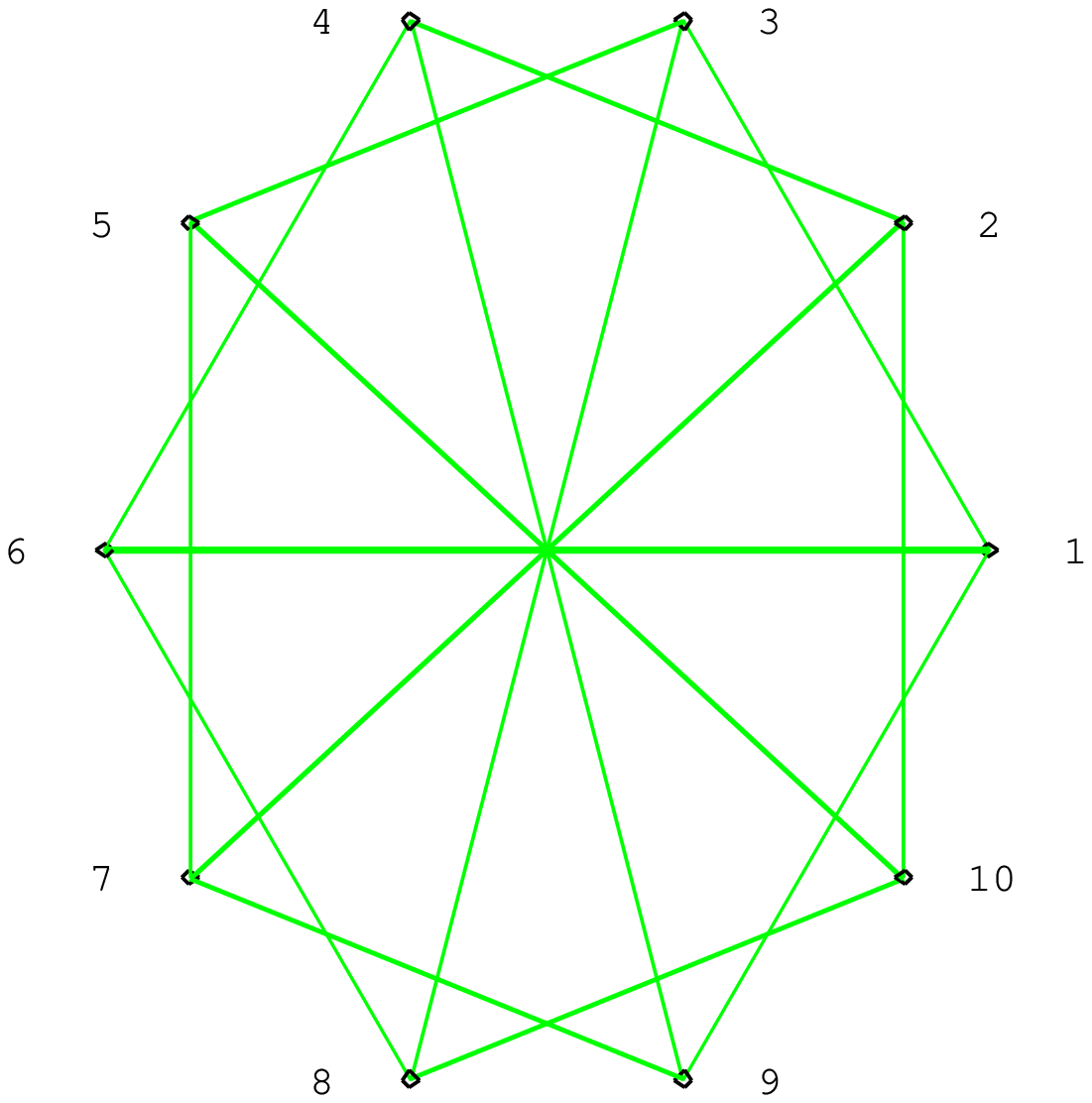} &
\includegraphics*[width=0.5\textwidth]
{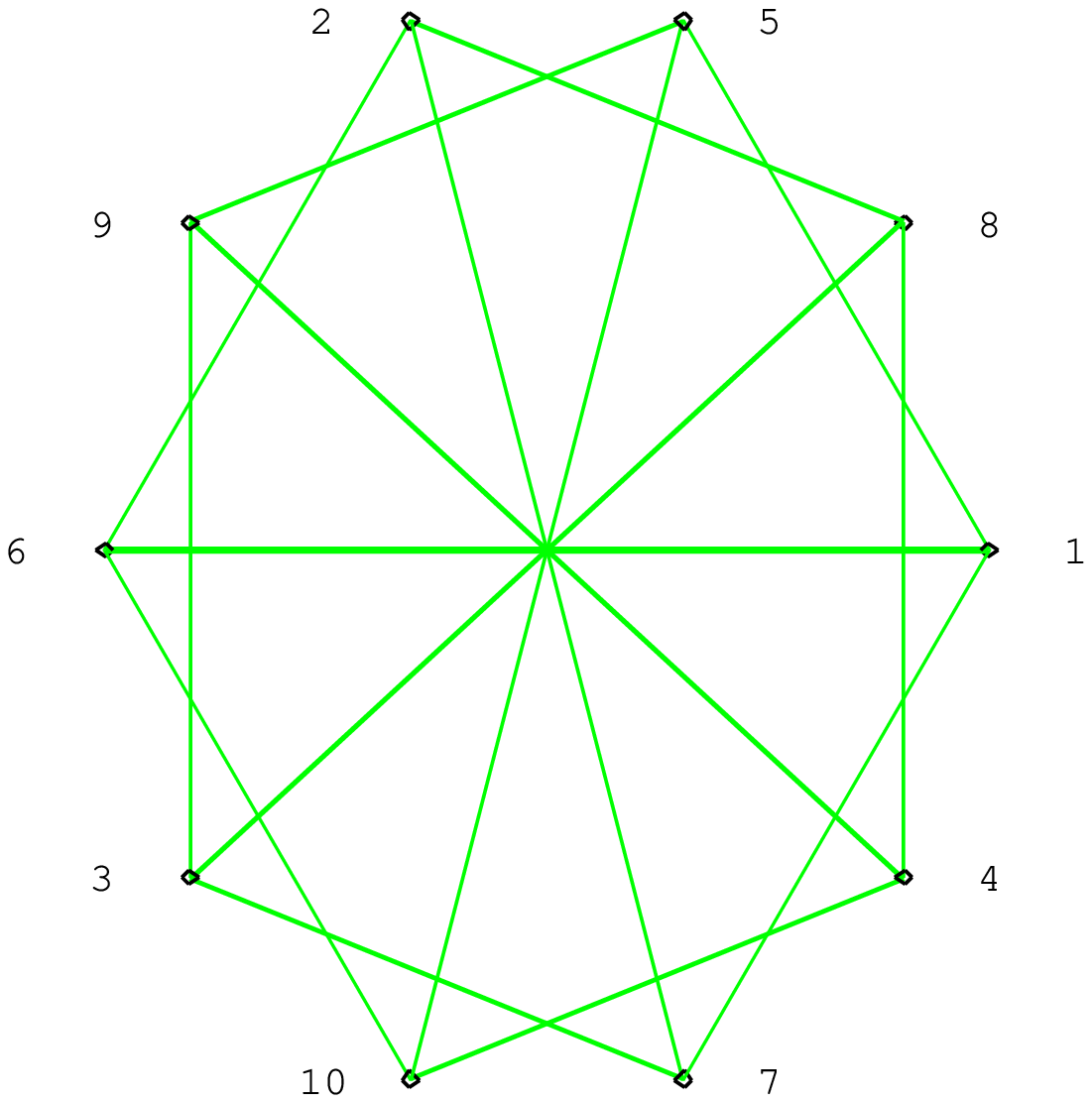}
\end{tabular}
$G_{\{2,5\}} \simeq G_{\{4,5\}}$ and 
$G_{\{1,2,3\}} \simeq G_{\{1,4,5\}}$.
\end{center}

\subsubsection{Inclusion poset}
The graph properties we study are monotone,
so it is important to know which graphs are
subgraphs of others. If $D\subset D'
\subseteq \{1,2,3,4,5\}$, then obviously
$G_D<G_{D'}$, but otherwise?

\begin{center}
\begin{tabular}{cc}
\includegraphics*[width=0.5\textwidth]
{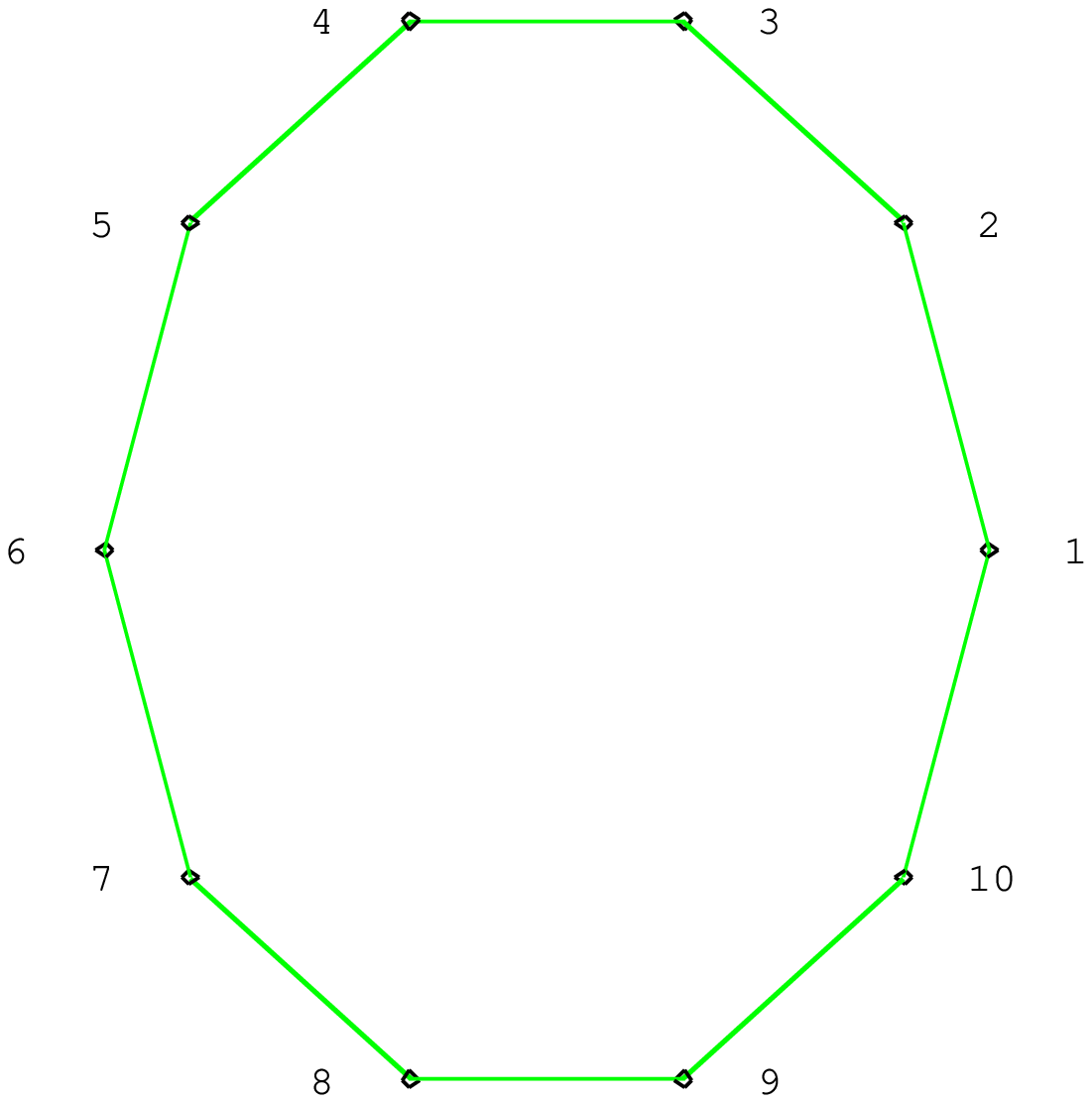} &
\includegraphics*[width=0.5\textwidth]
{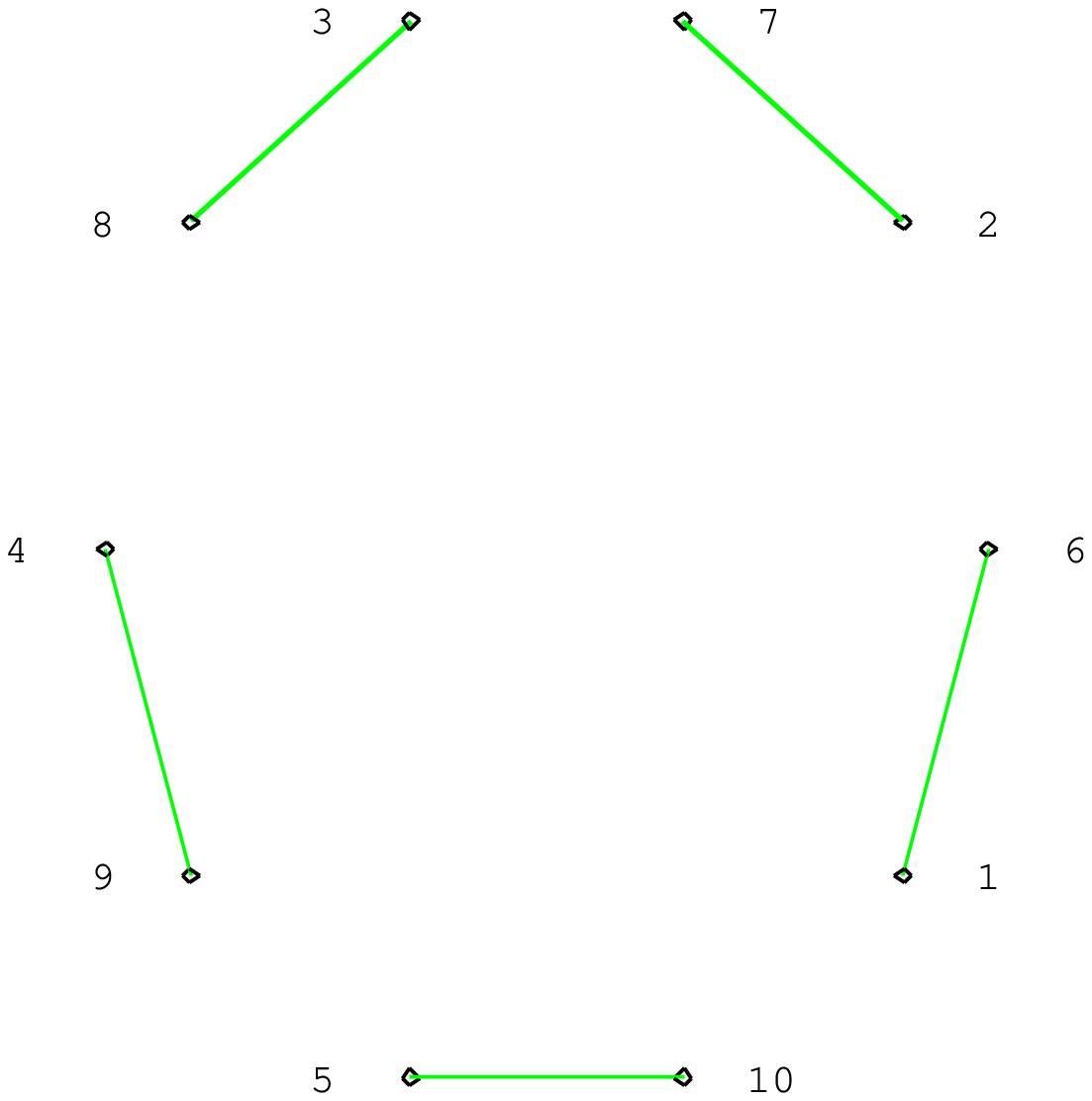}
\end{tabular}
$G_{\{1\}} \simeq G_{\{3\}}>G_{\{5\}}$ and 
$G_{\{1,2,3,4\}}>
G_{\{2,3,4,5\}} \simeq G_{\{1,2,4,5\}}$.
\end{center}

\begin{center}
\begin{tabular}{cc}
\includegraphics*[width=0.5\textwidth]
{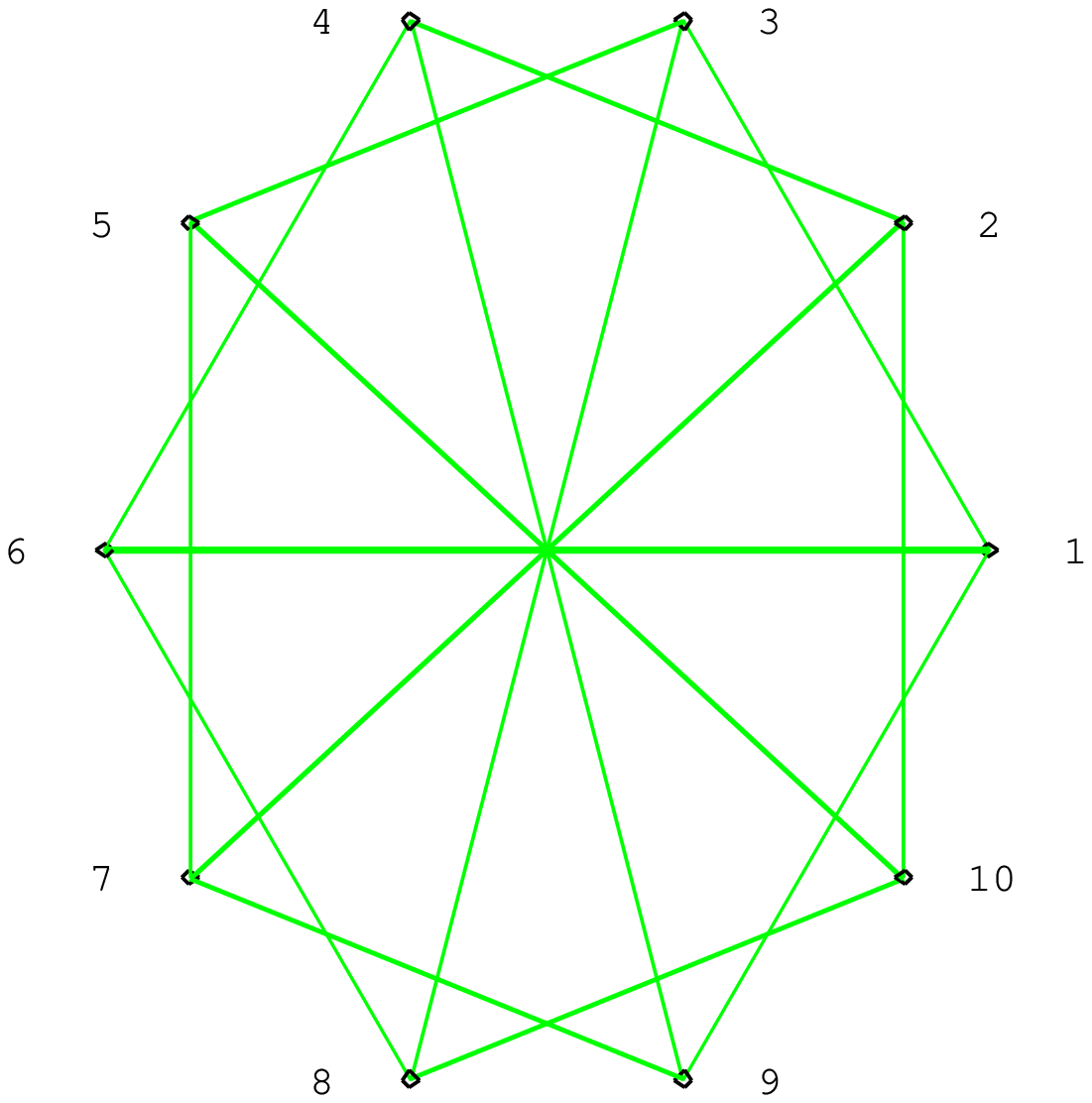} &
\includegraphics*[width=0.5\textwidth]
{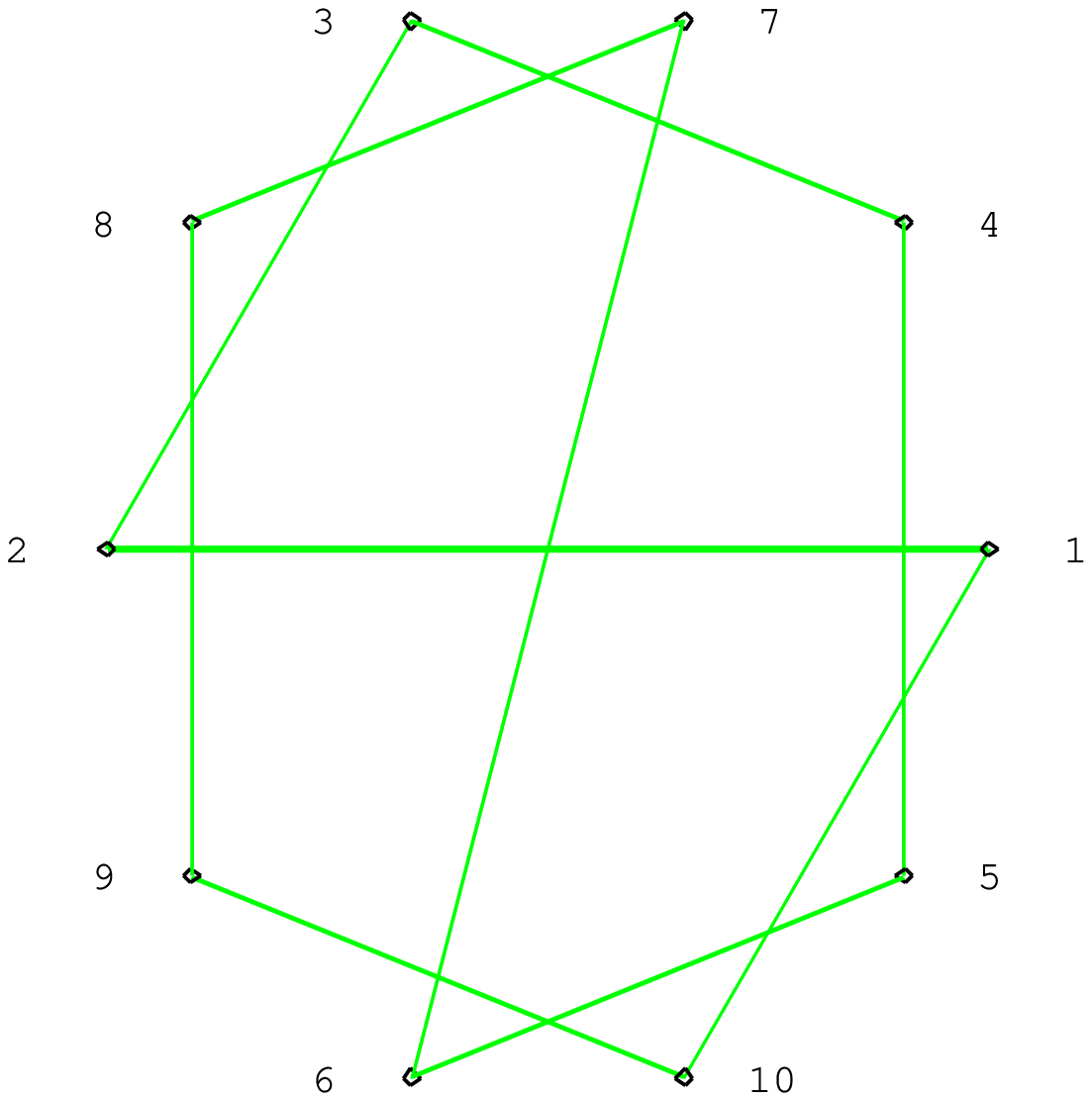}
\end{tabular}
$G_{\{2,5\}} \simeq G_{\{4,5\}}>
G_{\{1\}} \simeq G_{\{3\}}$ and
$G_{\{2,3,4,5\}} \simeq G_{\{1,2,4,5\}}>
G_{\{1,2,3\}} \simeq G_{\{1,3,4\}}$.
\end{center}

\begin{center}
\begin{tabular}{cc}
\includegraphics*[width=0.5\textwidth]
{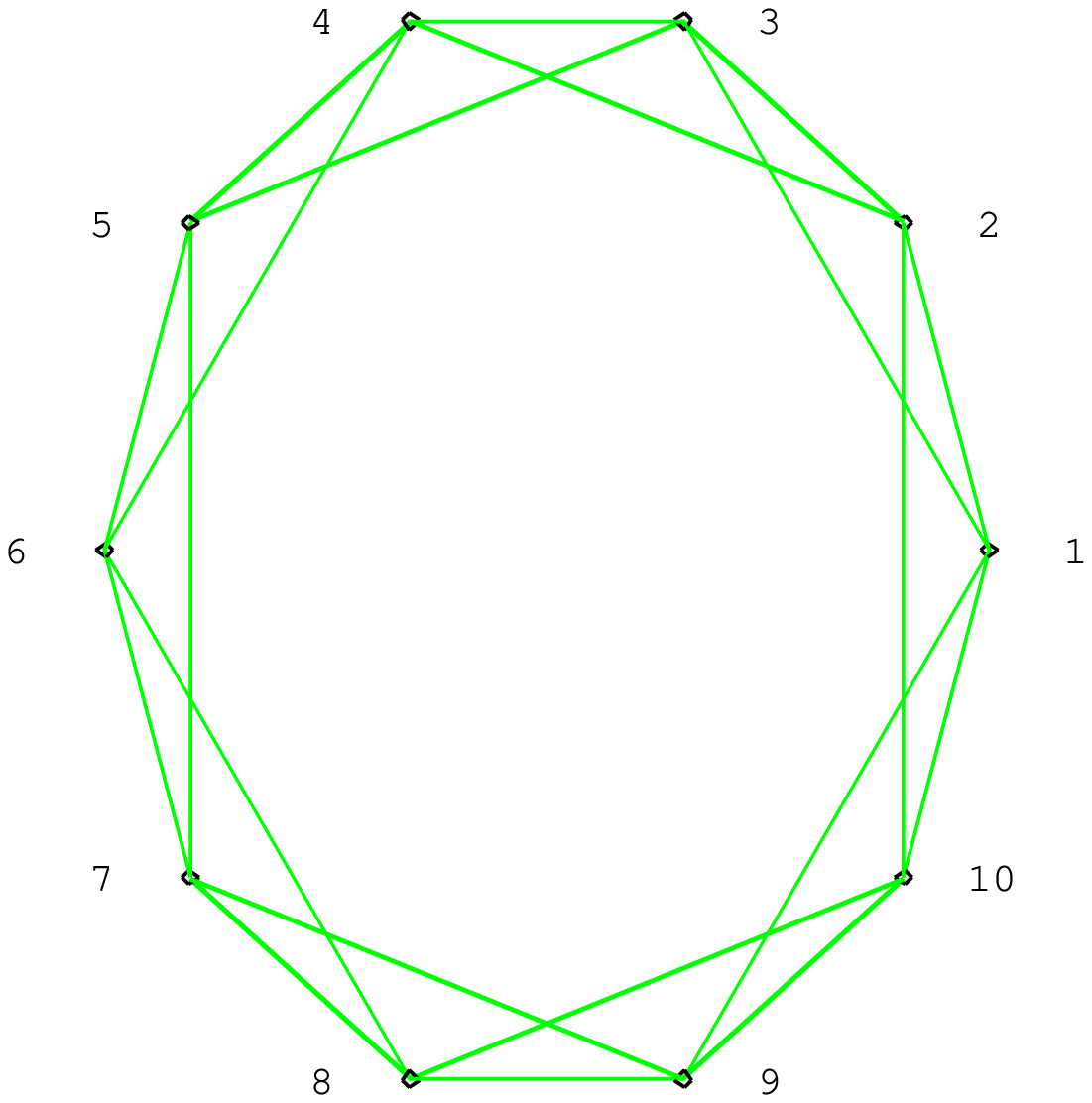} &
\includegraphics*[width=0.5\textwidth]
{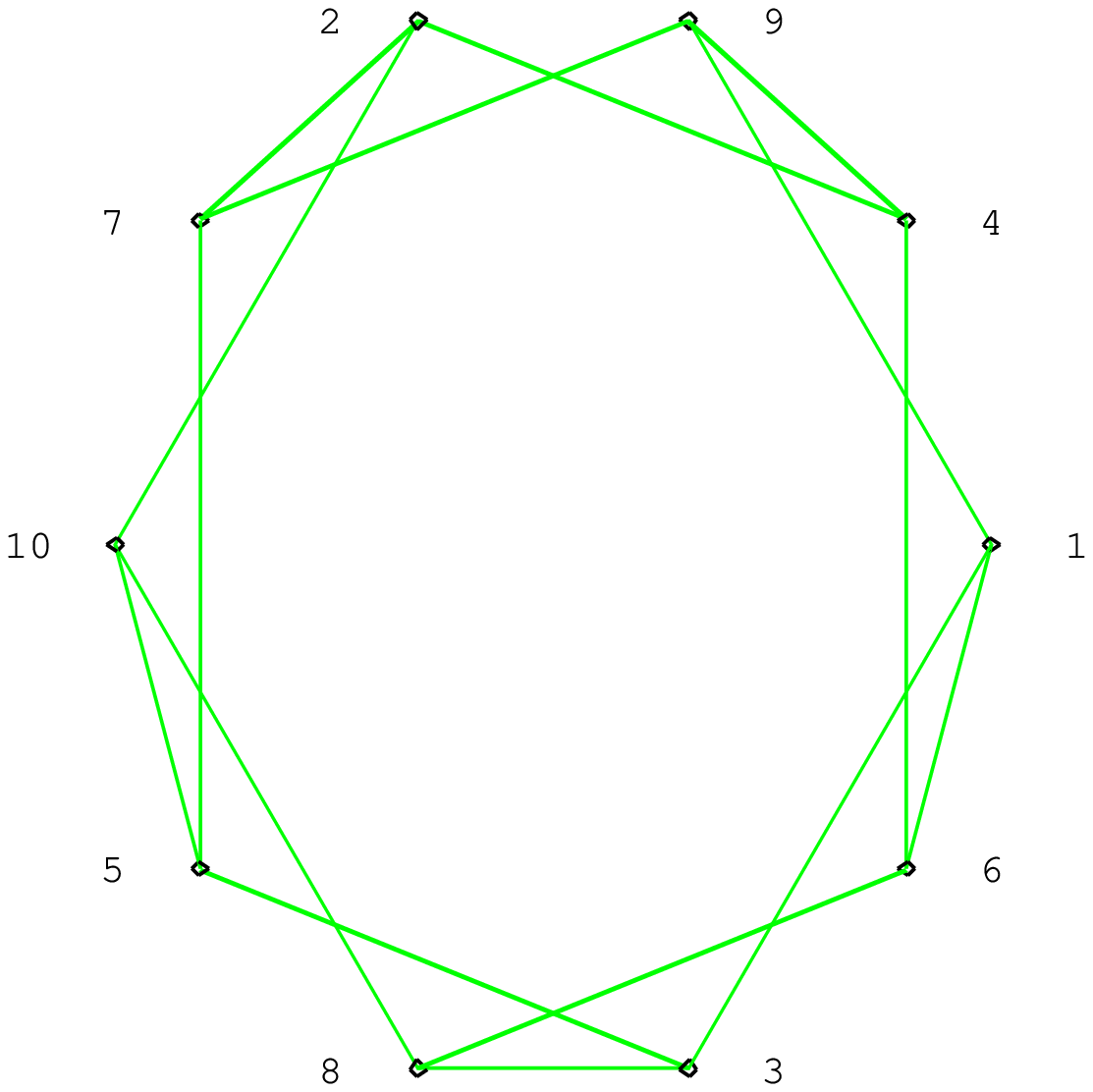}
\end{tabular}
$G_{\{1,2\}} \simeq G_{\{3,4\}}>
G_{\{2,5\}} \simeq G_{\{4,5\}}$ and
$G_{\{1,2,3\}} \simeq G_{\{1,3,4\}}>
G_{\{1,2,5\}} \simeq G_{\{3,4,5\}}$.
\end{center}

\begin{center}
\begin{tabular}{cc}
\includegraphics*[width=0.5\textwidth]
{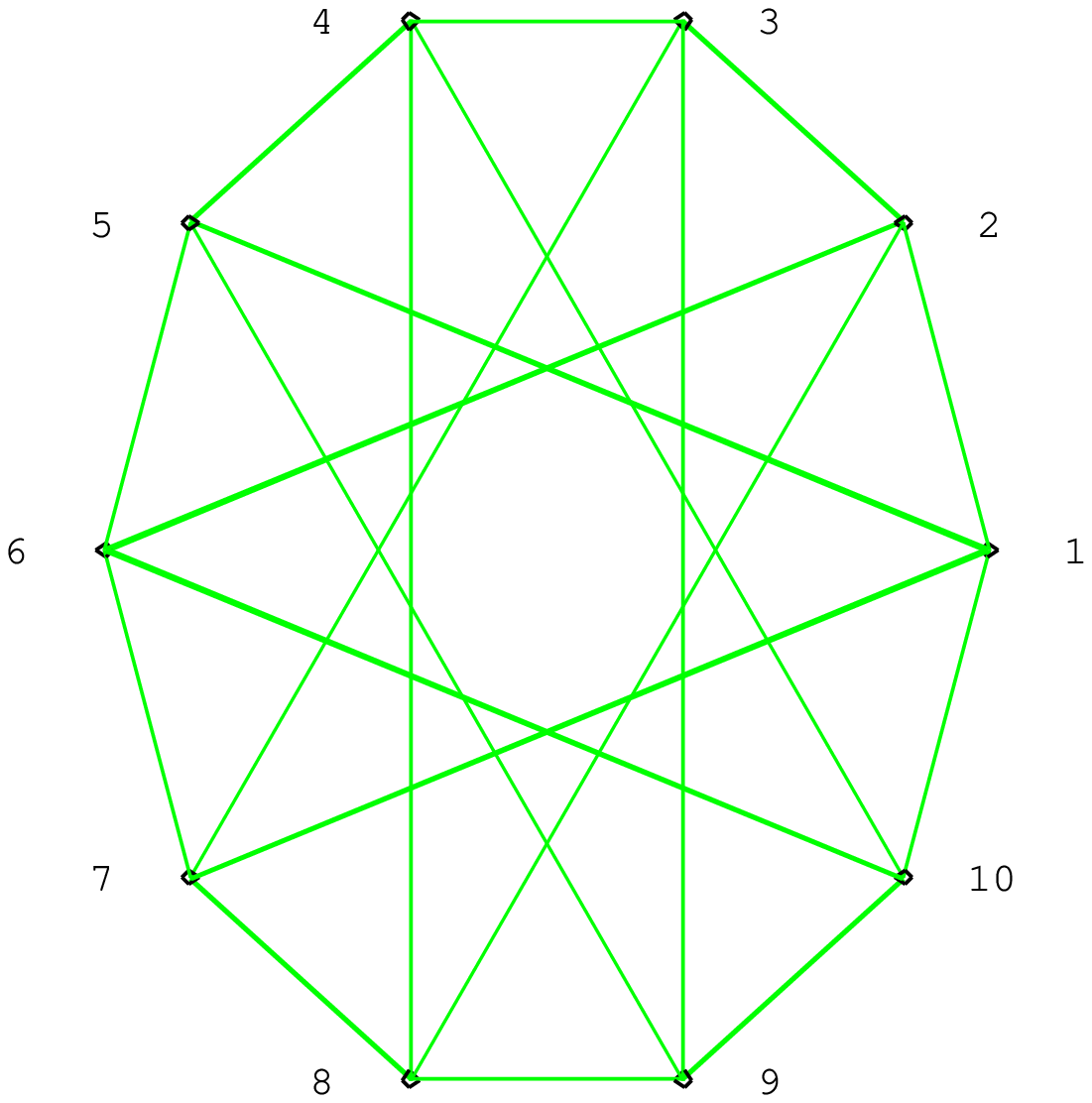} &
\includegraphics*[width=0.5\textwidth]
{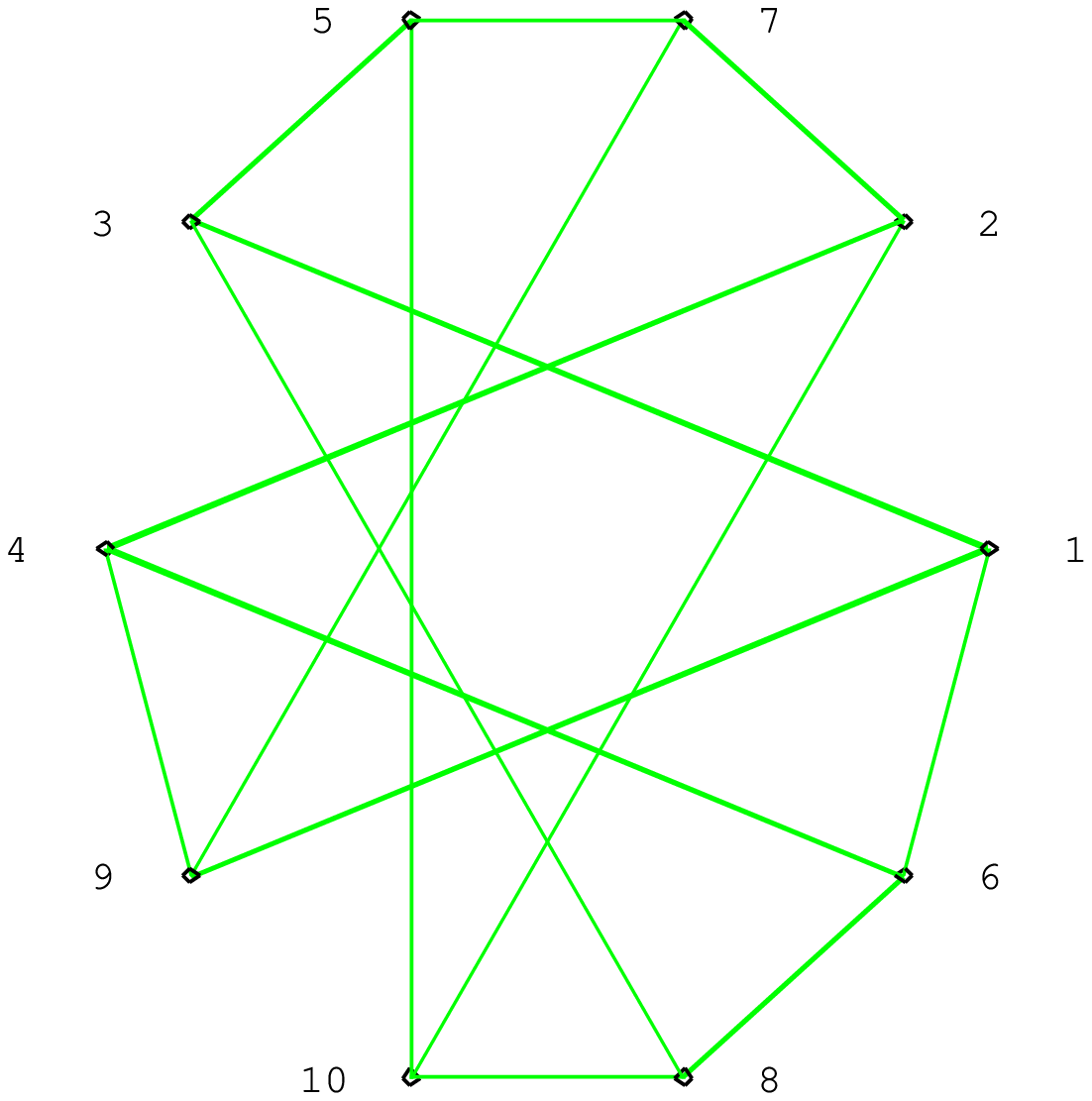}
\end{tabular}
$G_{\{1,4\}} \simeq G_{\{2,3\}}>
G_{\{2,5\}} \simeq G_{\{4,5\}}$ and
$G_{\{1,2,3\}} \simeq G_{\{1,3,4\}}>
G_{\{1,4,5\}} \simeq G_{\{2,3,5\}}$.
\end{center}

\begin{center}
\begin{tabular}{cc}
\includegraphics*[width=0.5\textwidth]
{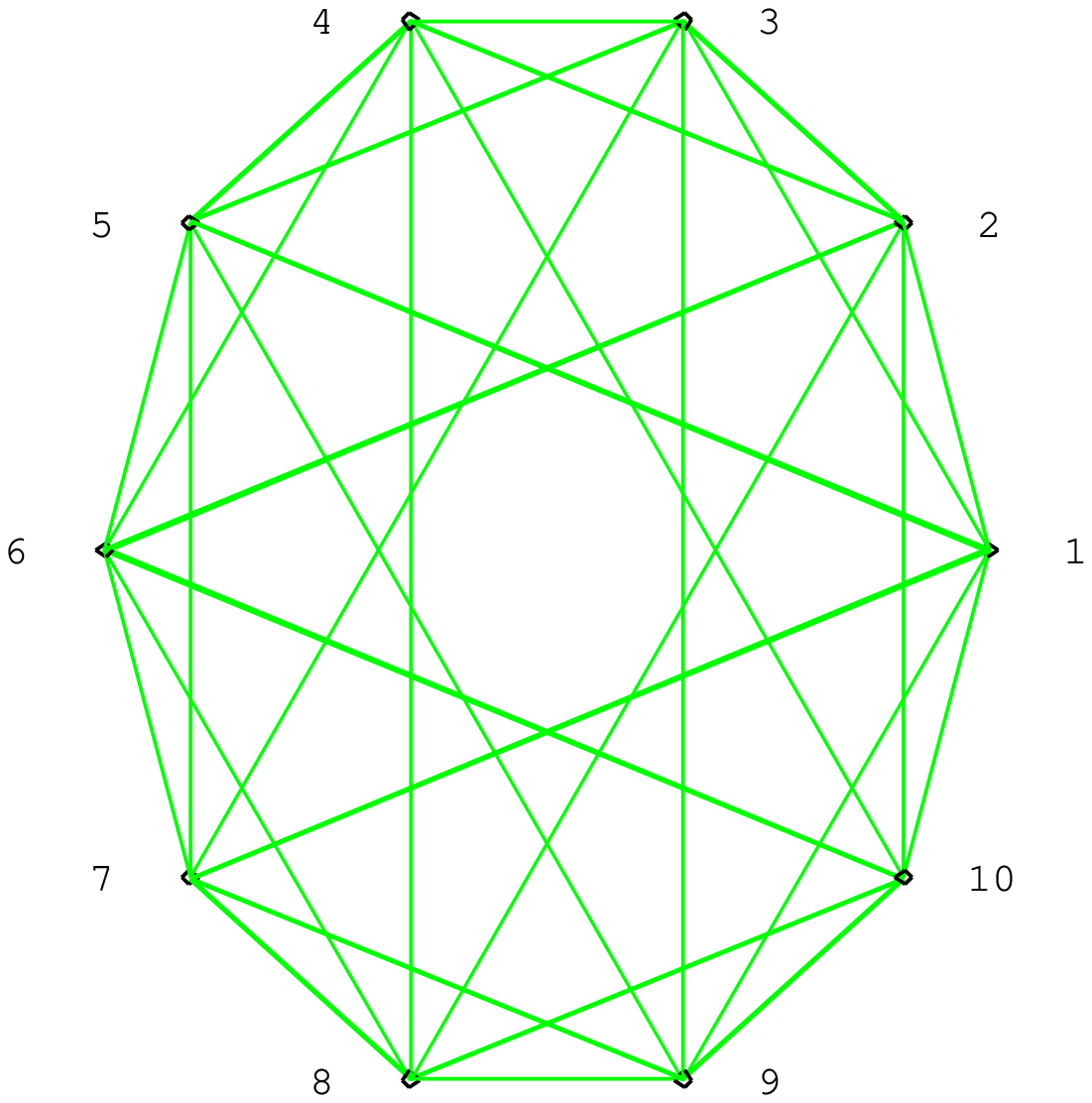} &
\includegraphics*[width=0.5\textwidth]
{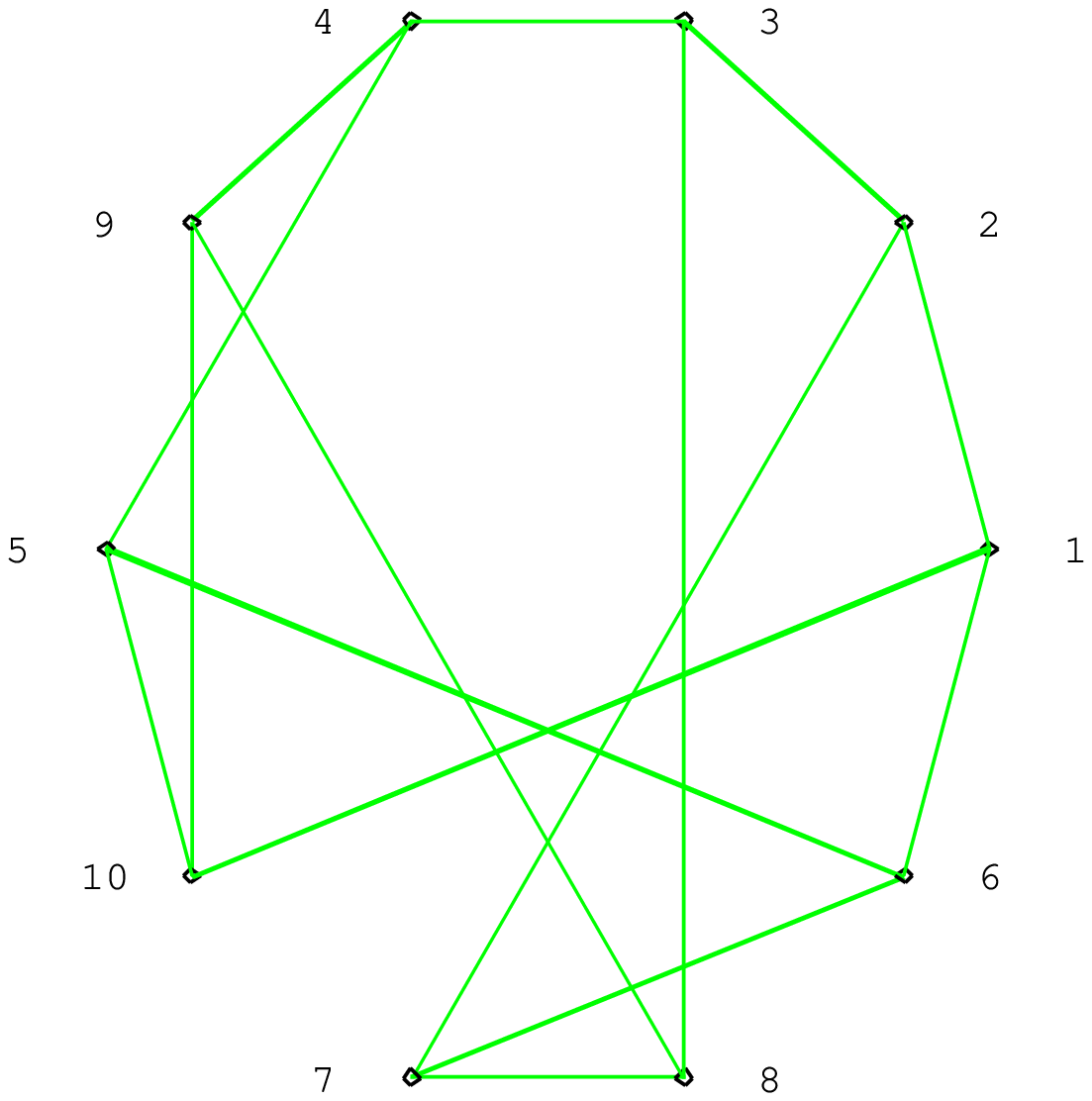}
\end{tabular}
$G_{\{1,2,4\}} \simeq G_{\{2,3,4\}}>
G_{\{1,5\}} \simeq G_{\{3,5\}}$
\end{center}

\begin{center}
\begin{tabular}{cc}
\includegraphics*[width=0.5\textwidth]
{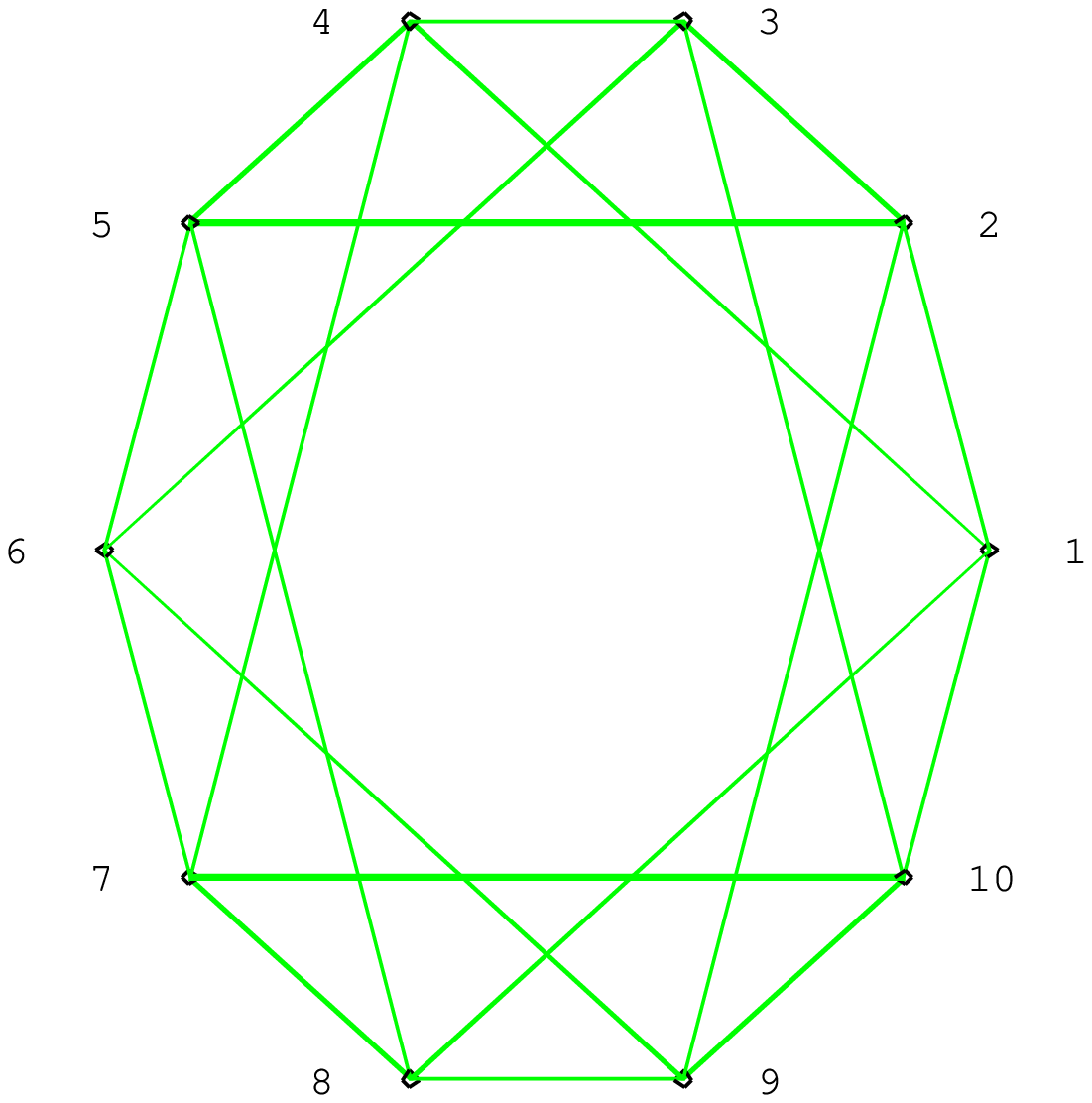} &
\includegraphics*[width=0.5\textwidth]
{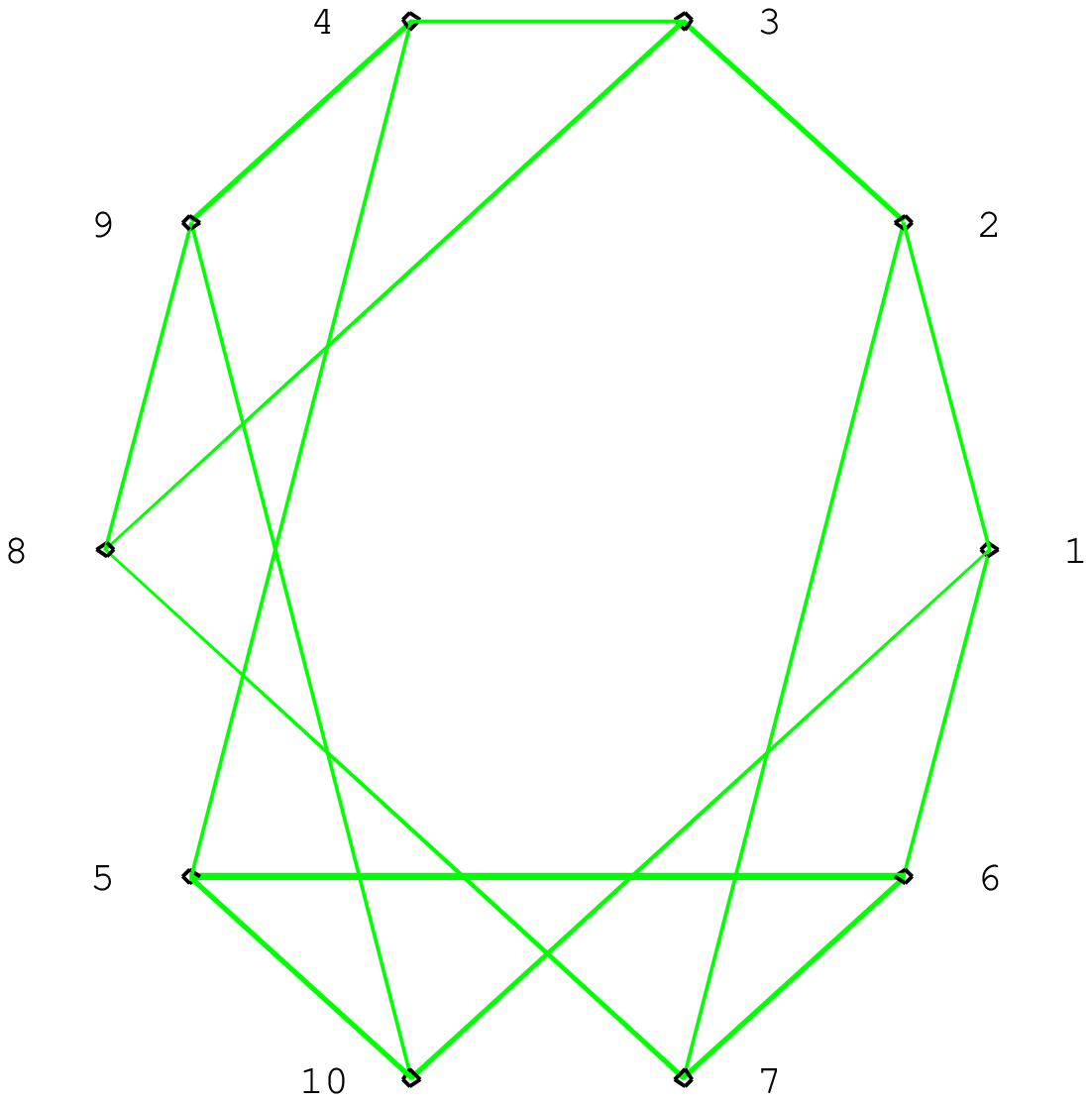}
\end{tabular}
$G_{\{1,3\}}>
G_{\{1,5\}} \simeq G_{\{3,5\}}$ and
$G_{\{1,2,4\}} \simeq G_{\{2,3,4\}}>
G_{\{2,4,5\}}$.
\end{center}

\begin{center}
\begin{tabular}{cc}
\includegraphics*[width=0.5\textwidth]
{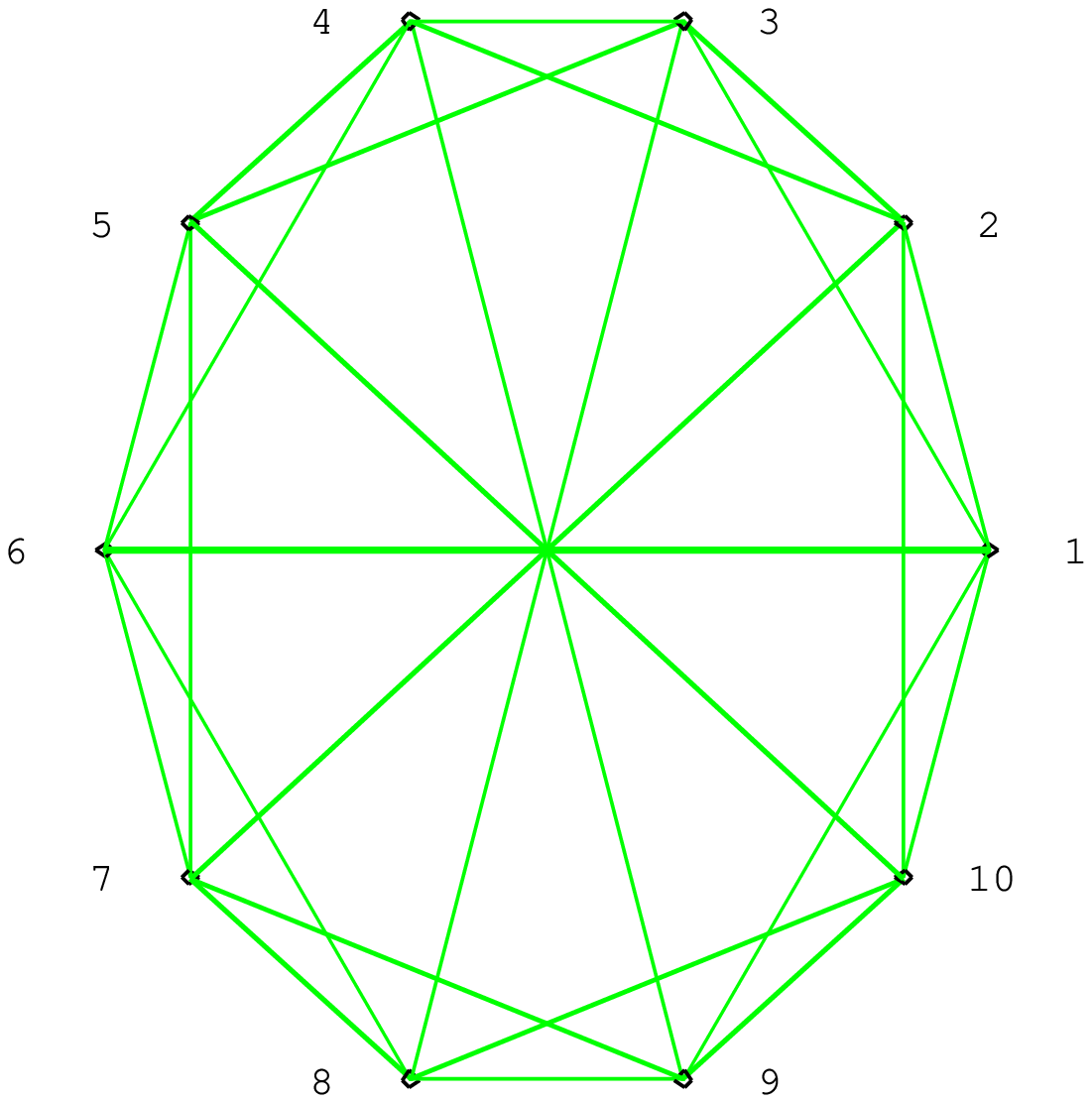} &
\includegraphics*[width=0.5\textwidth]
{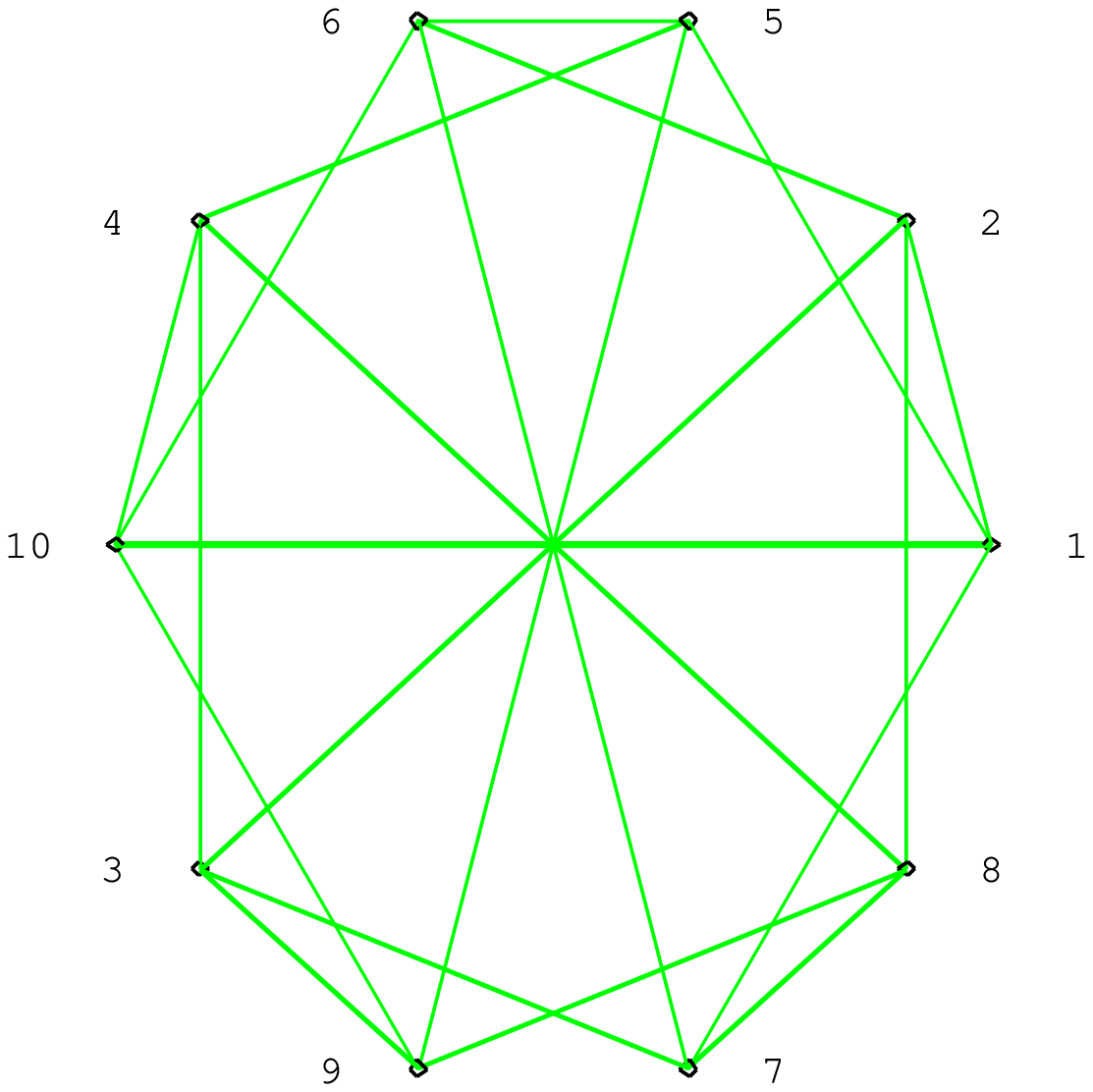}
\end{tabular}
$G_{\{1,2,5\}} \simeq G_{\{3,4,5\}}>
G_{\{1,4\}} \simeq G_{\{2,3\}}$ and
$G_{\{1,4,5\}} \simeq G_{\{2,3,5\}}>
G_{\{1,2\}} \simeq G_{\{3,4\}}$.
\end{center}

The inclusion relation give a poset which
is drawn in the next figure. Grey lines
are trivial inclusions, and black lines are
inclusions listed above.
\begin{center}
\includegraphics*{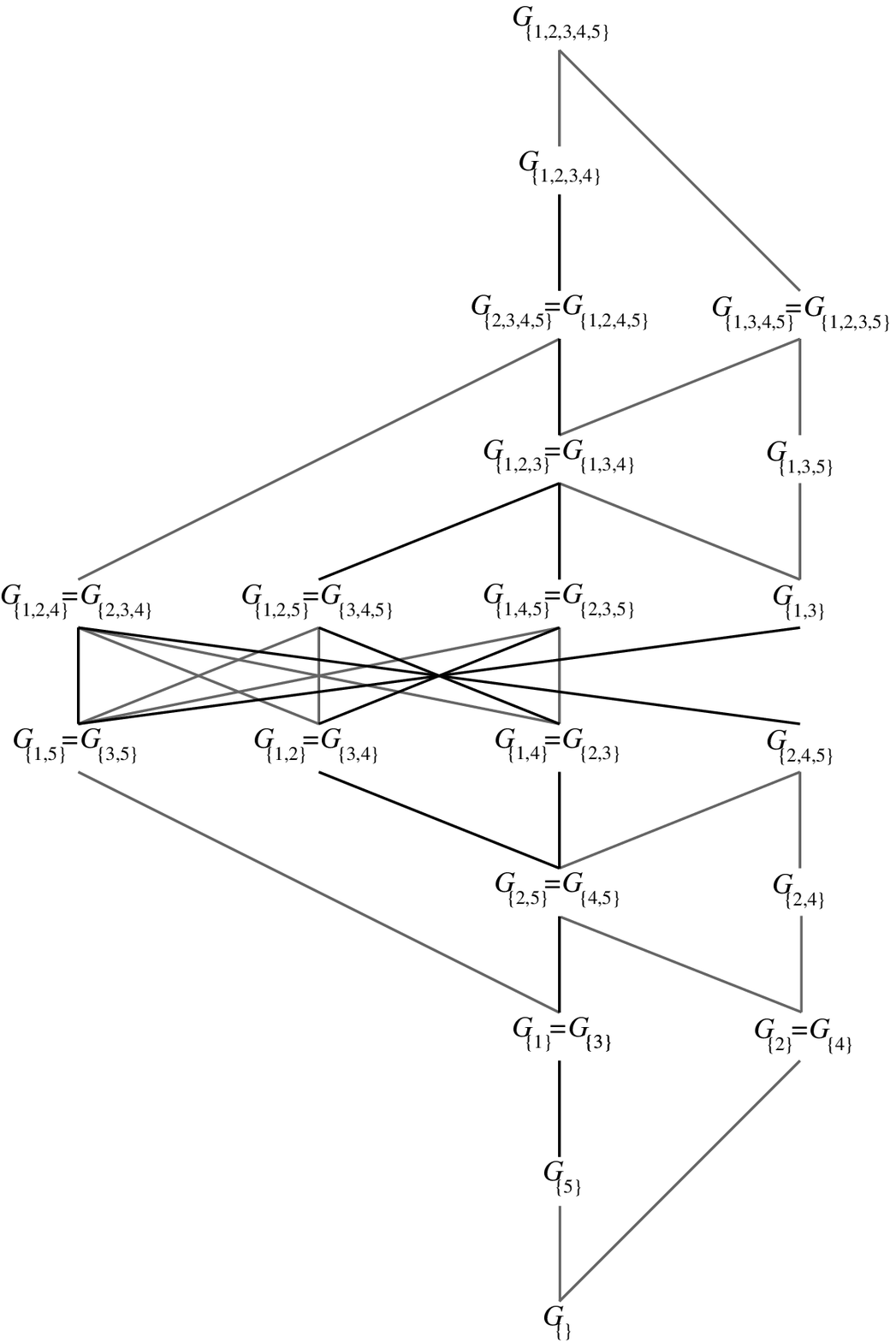}
\end{center}

\subsection{The Petersen graph}
The non Cayley transitive graphs are
the Petersen graph and its complement.

\begin{center}
\begin{tabular}{cc}
\includegraphics*[width=0.5\textwidth]
{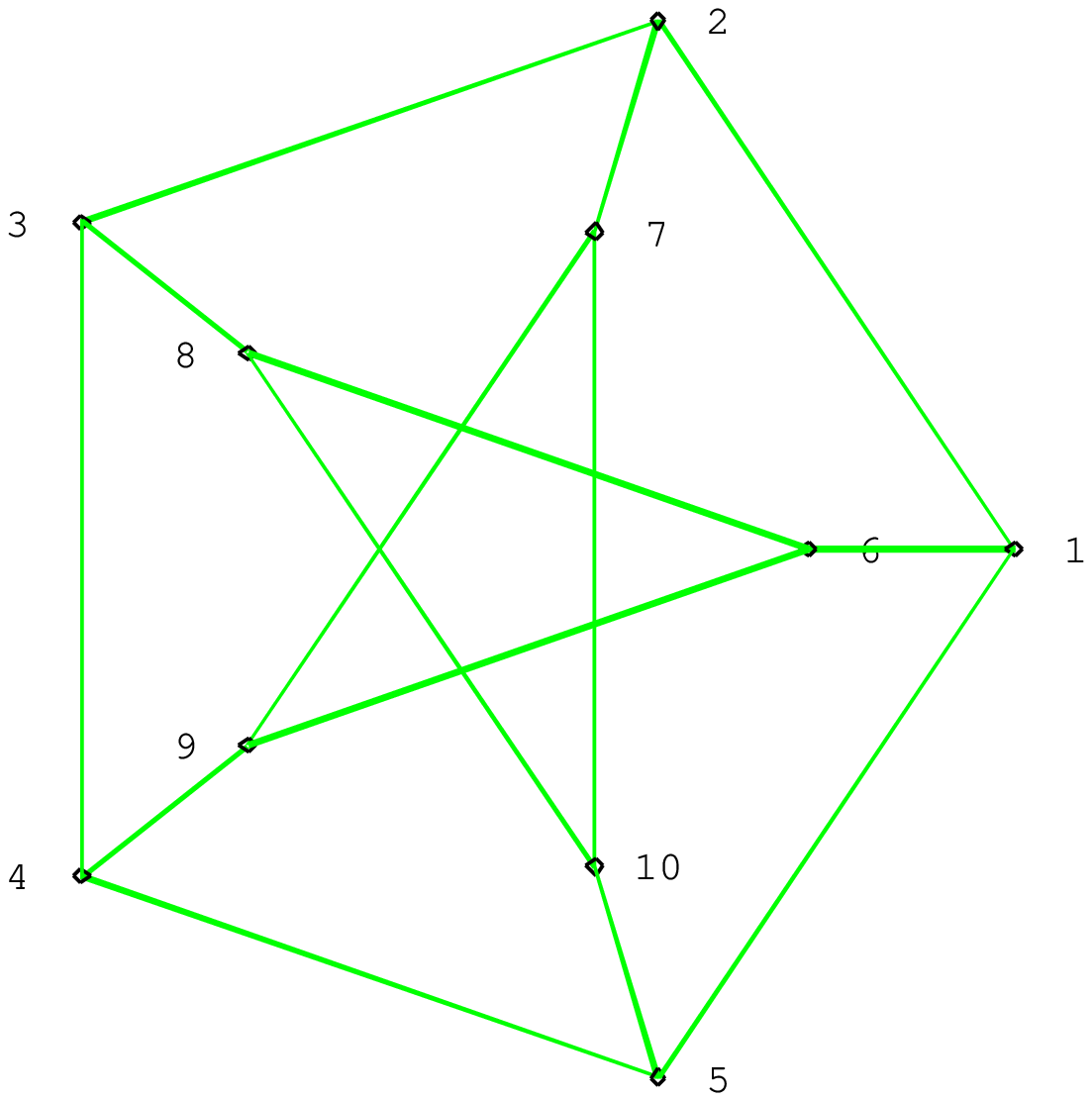} &
\includegraphics*[width=0.5\textwidth]
{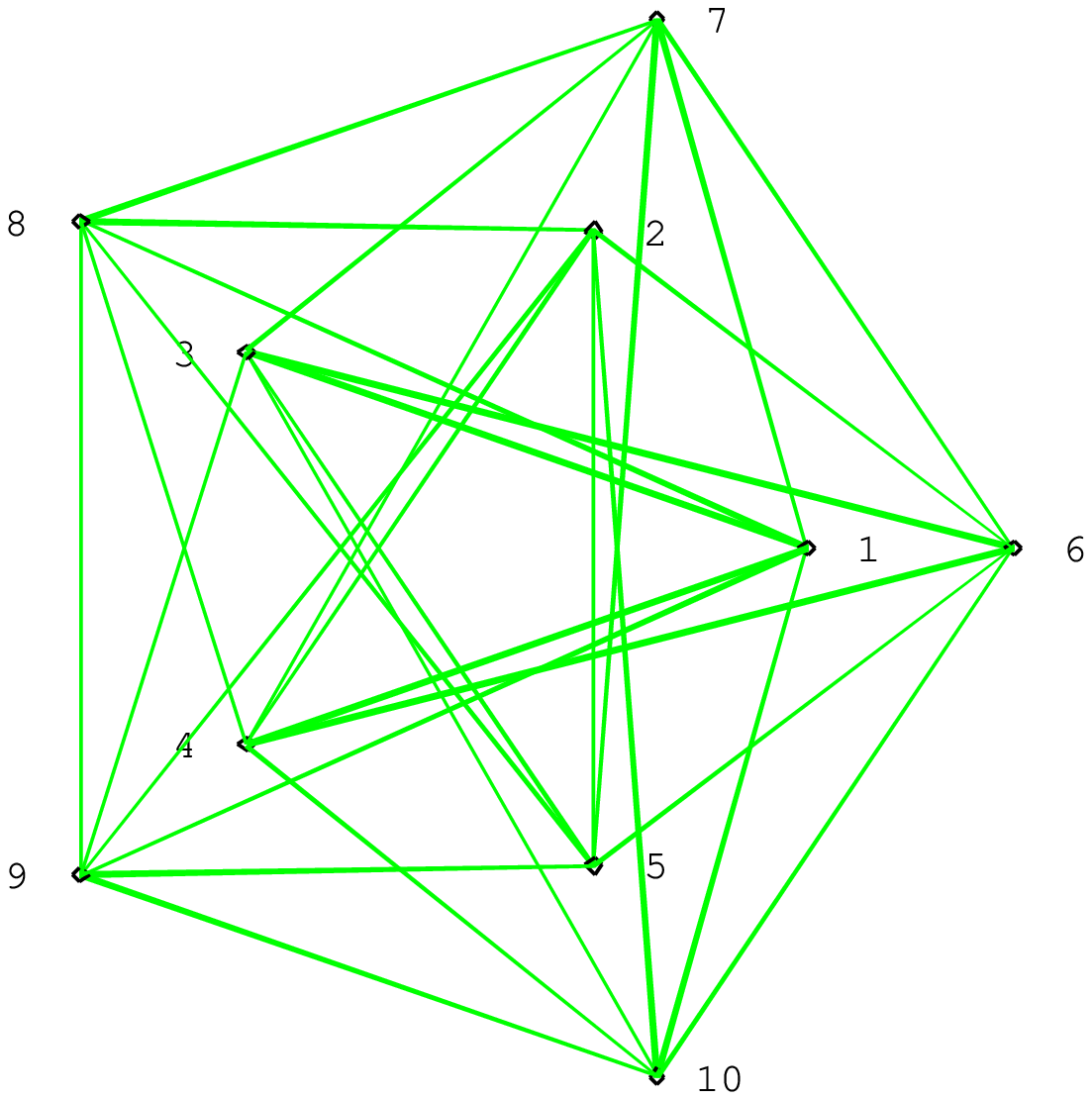}\\
The Petersen graph $P$.&
The complement of the
Petersen graph, $\bar{P}.$
\end{tabular}
\end{center}
From the figure above it is clear that $P<
\bar{P}$. The other relations are:
\begin{center}
\begin{tabular}{cc}
\includegraphics*[width=0.5\textwidth]
{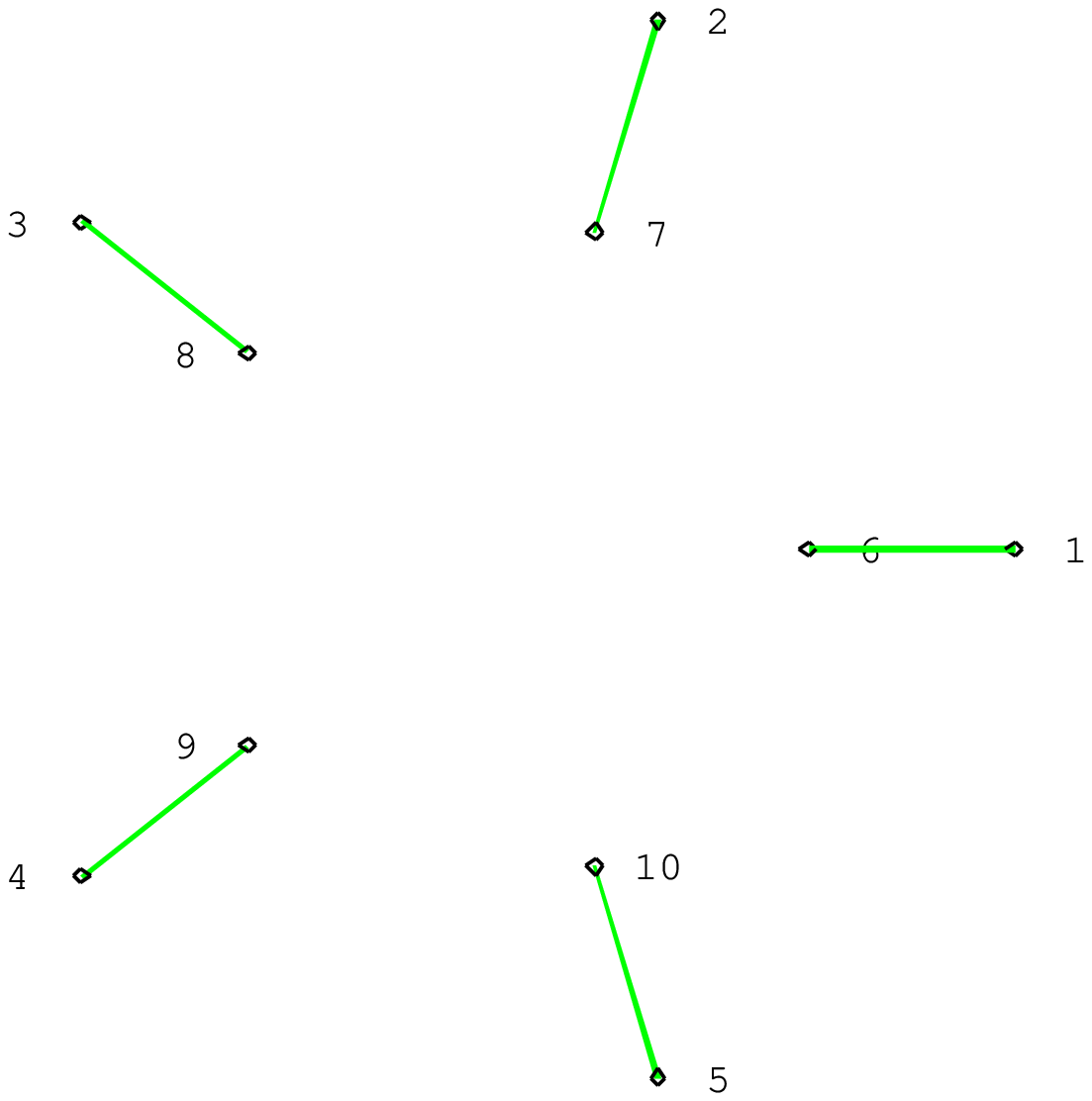} &
\includegraphics*[width=0.5\textwidth]
{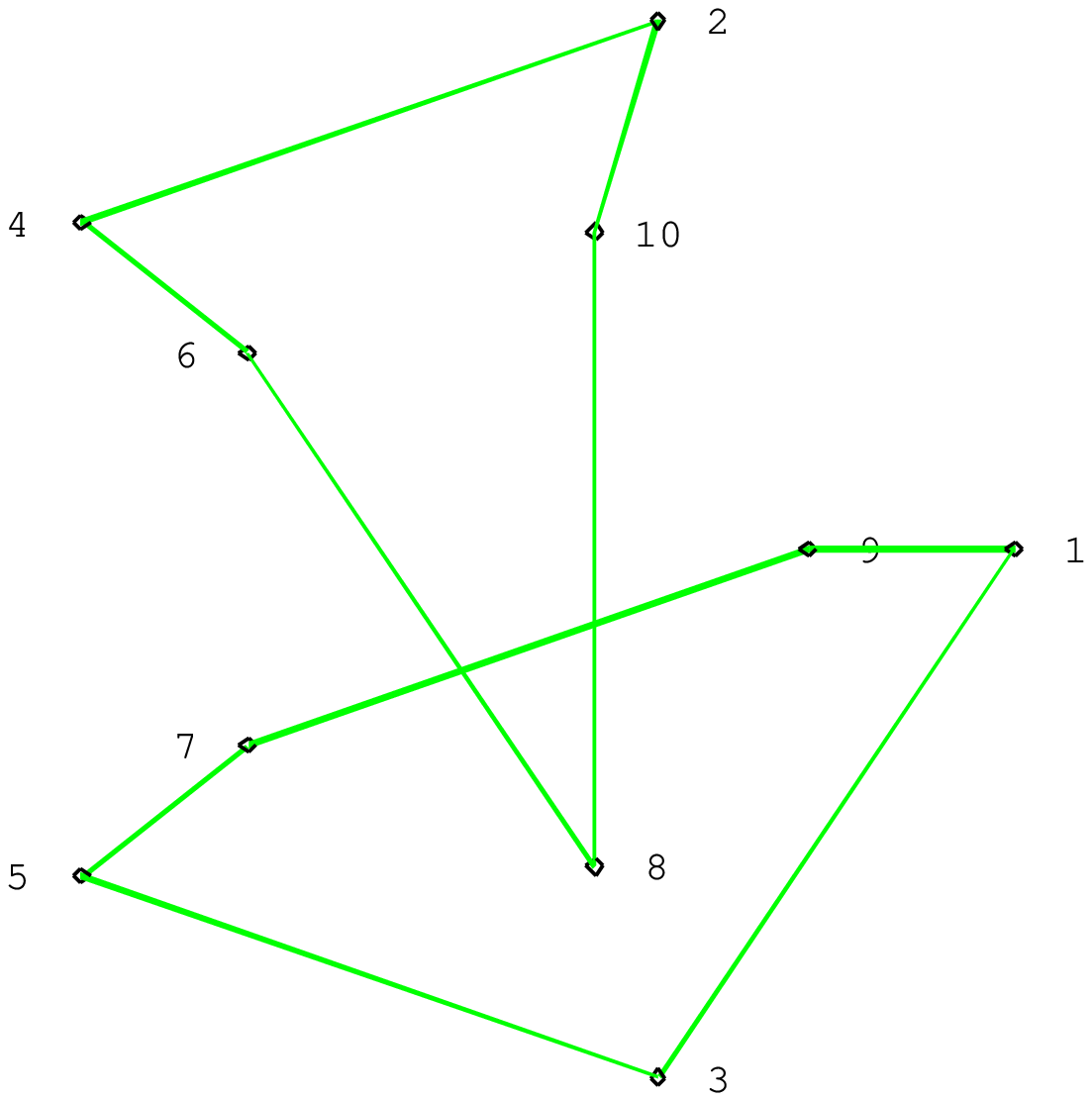}\\
$G_{\{5\}}<P$ &
$G_{\{2\}}\simeq G_{\{4\}}<P$ \\
$\bar{P}<G_{\{1,2,3,4\}}$ &
$\bar{P}<G_{\{1,2,3,5\}}\simeq G_{\{1,3,4,5\}}$\\
\end{tabular}
\end{center}

\begin{center}
\begin{tabular}{cc}
\includegraphics*[width=0.5\textwidth]
{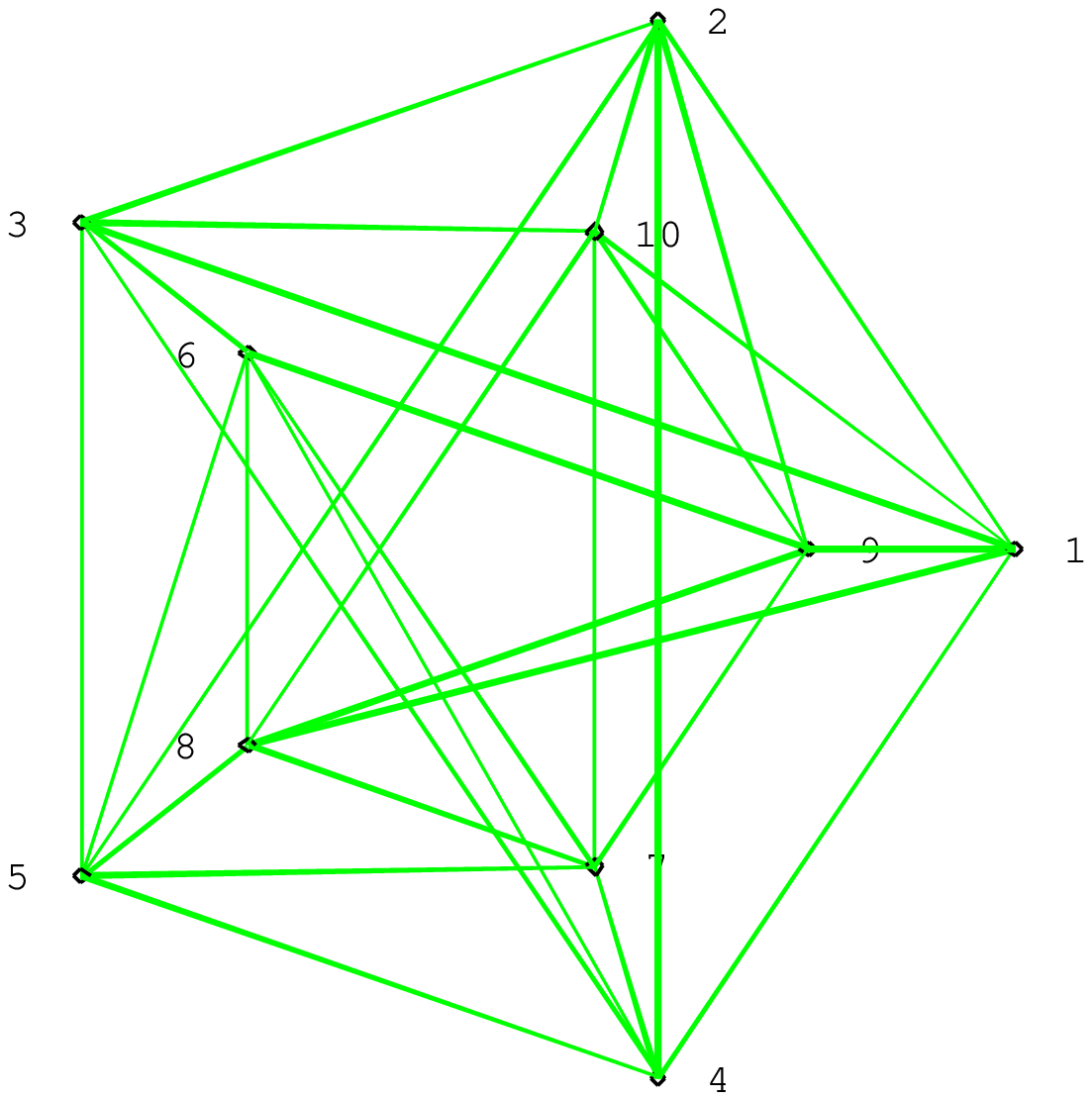} &
\includegraphics*[width=0.5\textwidth]
{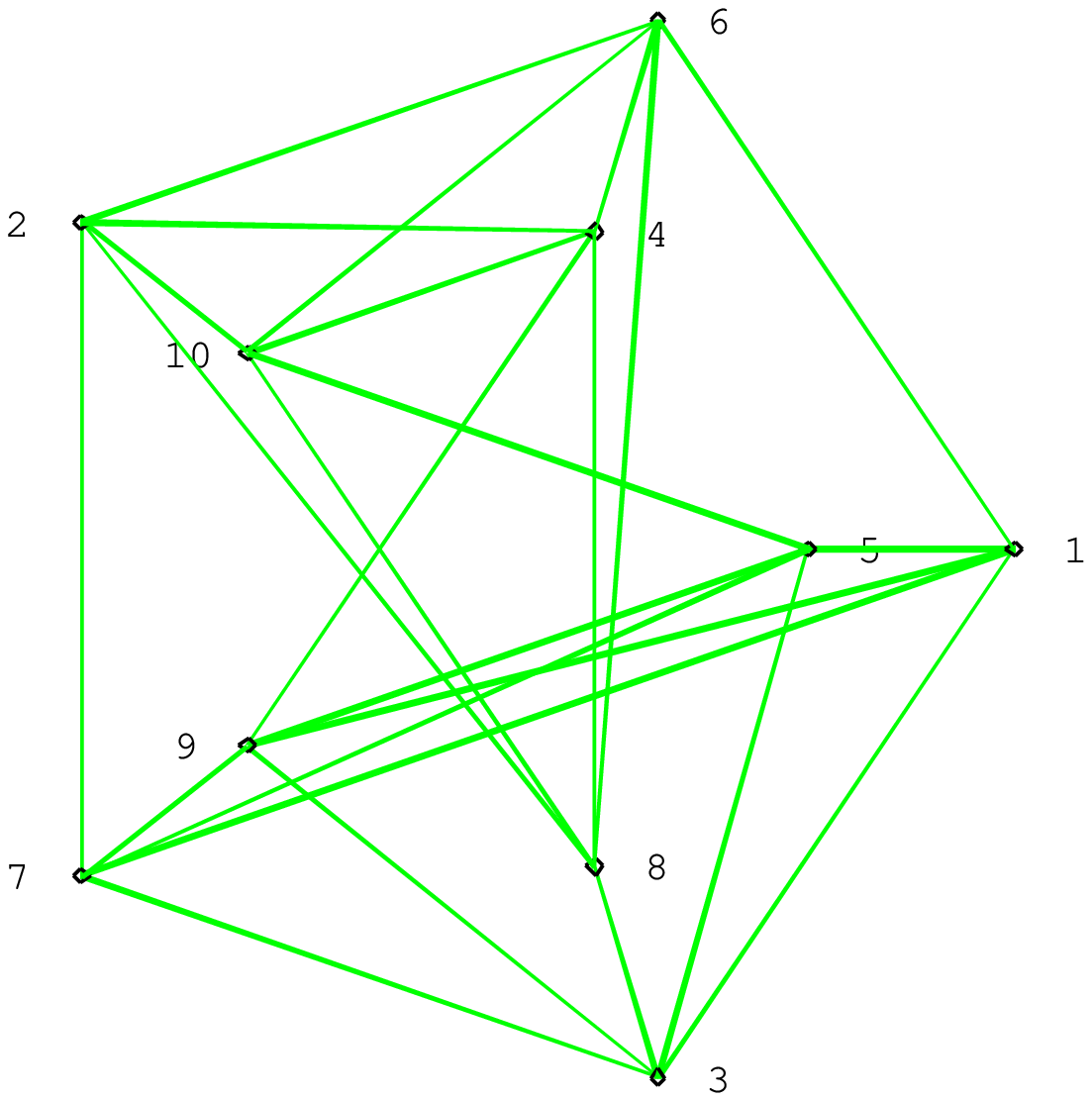}\\
$P<G_{\{1,2,3\}}\simeq G_{\{1,3,4\}}$ &
$P<G_{\{2,4,5\}} $ \\
$G_{\{2,5\}}\simeq G_{\{4,5\}}<\bar{P}$ &
$G_{\{1,3\}}<\bar{P}$\\
\end{tabular}
\end{center}

Inserting $P$ and $\bar{P}$ in the inclusion
poset of the Cayley graphs would be a mess.
But selecting the elements with a cover relation
to $P$ and $\bar{P}$ gives the poset:
\begin{center}
\includegraphics*{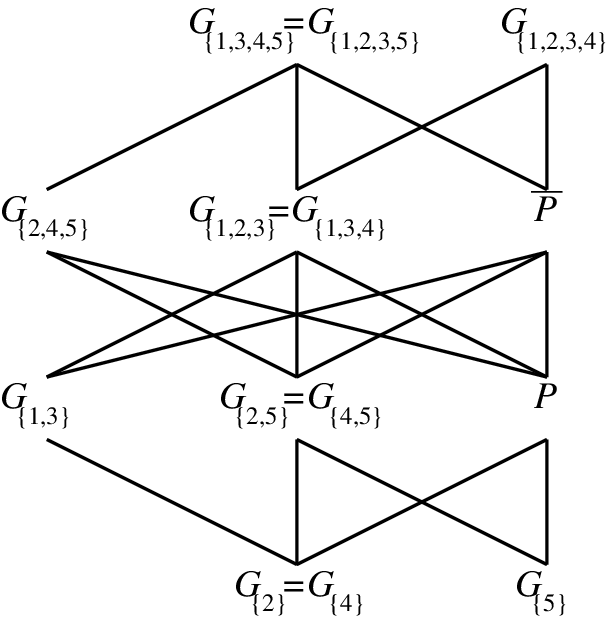}
\end{center}

\section{Using the topological method}\label{se}
If $G$ is a graph on 10 vertices,
and $\Delta$ a simplicial complex, then
define the indicator $i_G$ as: $i_G=1$
if $G\in \Delta$ and $i_G=0$ if $G
\not \in \Delta$. To simplify notation
for the Cayley graphs, we write
$i_{124}$ instead of $i_{G_{\{1,2,4\}}}$
for example. Note that by monotonicity
$i_G\leq i_{G'}$ if $G\supseteq G'$. 
The two posets in the last section carries
over to the indicators, and that is how
we will use them.
 
Let $\mathfrak{C}$ be the set of
nonevasive nontrivial graph complexes on 
10 vertices. Note that $i_{12345}=0$ for
all nontrivial complexes.



\begin{lemma}\label{lemma:T24}
For any simplicial complex in $\mathfrak{C}$,
$i_5=1$ and $i_{1234}=0$.
\end{lemma}
\begin{proof}
Let
\begin{itemize}
\item
 $\Gamma=\langle
(1\: 3\: 5\: 7\: 9)(2\: 4\: 6\: 8\: 10),
(1\: 7\: 9\: 3)(2\: 4\: 8\: 6),
(2\: 7)(5\: 10)
\rangle$.
\item Define a homomorphism $\phi$ from $\Gamma$
onto $\mathbb{Z}_4$ by $\phi((1\: 3\: 5\: 7\: 9)
(2\: 4\: 6\: 8\: 10))=0$, $\phi((1\: 7\: 9\: 3)
(2\: 4\: 8\: 6))=1$ and $\phi((2\: 7)(5\: 10))
=0$.
\item
$\Gamma'= \mathrm{Ker}(\phi)=
\langle
(1\: 3\: 5\: 7\: 9)(2\: 4\: 6\: 8\: 10),
(2\: 7)(5\: 10)
\rangle$.
\item Define a homomorphism $\phi'$ from 
$\Gamma'$ onto $\mathbb{Z}_5$ by
$\phi'((1\: 3\: 5\: 7\: 9)(2\: 4\: 6\: 
8\: 10))=3$, and $\phi'((2\: 7)(5\: 10)=0$.
\item
$\Gamma''= \mathrm{Ker}(\phi')=
\langle
(2\: 7)(5\: 10),(4\: 9)(5\: 10),(3\: 8)(5\: 10),
(1\: 6)(2\: 7)(3\: 8)(5\: 10)\rangle$.
\end{itemize}
Then
\begin{itemize}
\item $\Gamma''\triangleleft
         \Gamma'\triangleleft \Gamma$,
\item $|\Gamma''|=2^4$,
\item $\Gamma'/\Gamma''=\mathbb{Z}_5$ is cyclic,
\item $|\Gamma'/\Gamma''|=|\mathbb{Z}_4|=2^2
\Rightarrow q=2$.
\end{itemize}
By theorem \ref{theorem:2}, $\chi(
\Delta^\Gamma)\equiv 1\: \mathrm{mod}\: q$.
The vertices of $\Delta^\Gamma$ are 
$\{G_{\{1,2,3,4\}},G_{\{5\}}\}$, and
$G_{\{5\}}<G_{\{1,2,3,4\}}$, so $i_5=1$
and $i_{1234}=0$ is the only way to achive
$\chi(\Delta^\Gamma)=i_5+i_{1234}\equiv
1\: \mathrm{mod}\: 2$.
\end{proof}


\begin{lemma}\label{lemma:T1}
For any simplicial complex in $\mathfrak{C}$,
\[
\begin{array}{l}
2i_{1}+2i_{2}-2i_{12}-i_{13}-
2i_{14}-2i_{15}-i_{24}-2i_{25}
+2i_{123}\\
+2i_{124}
+2i_{125}+i_{135}+2i_{145}
+ i_{245}-2i_{1235}-2i_{1245}=0
\end{array}
\]
\end{lemma}
\begin{proof}
Let $\Gamma=\langle (1\:2\:3\:4\:5\:6\:
7\:8\:9\:10)\rangle$ and $\Gamma'=\mathrm{Id}$.
By theorem \ref{theorem:1}, 
$\chi(\Delta^\Gamma)=1$. The vertex set of
$\Delta^\Gamma$ is $\{G_{\{1\}}, G_{\{2\}},
G_{\{3\}}, G_{\{4\}}, G_{\{5\}}\}$, so
\[
\begin{array}{rcl}
{\chi(\Delta^\Gamma)} & = &
i_1+i_2+i_3+i_4+i_5 \\
&& - i_{12}-i_{13}-i_{14}-i_{15}-i_{23}-
i_{24}-i_{25}-i_{34}-i_{35}-i_{45}\\
&& + i_{123}+i_{124}+i_{125}+i_{134}+i_{135}
+i_{145}+i_{234}+i_{235}+i_{245}+i_{345}\\
&& -i_{1234}-i_{1235}-i_{1245}-i_{1345}
-i_{2345}\\
& = & 2i_1+2i_2+1 \\
&& - 2i_{12}-i_{13}-2i_{14}-2i_{15}-
i_{24}-2i_{25}\\
&& + 2i_{123}+2i_{124}+2i_{125}+i_{135}
+2i_{145}+i_{245}\\
&& -0-2i_{1235}-2i_{1245}\\
&=&1\\
\end{array}\]
\end{proof}


\begin{lemma}\label{lemma:T4}
For any simplicial complex in $\mathfrak{C}$,
\[
2i_2+i_{13}-i_{24}-2i_{25}-2i_{123}-i_{135}
+i_{245}+2i_{1235}=0
\]
\end{lemma}
\begin{proof}
Let
\begin{itemize}
\item $\Gamma=\langle
(1\: 3\: 5\: 7\: 9)(2\: 4\: 6\: 8\: 10),
(1\: 2\: 9\: 8)(3\: 6\: 7\: 4)(5\: 10)\rangle$. 
\item Define a homomorphism $\phi$ from $\Gamma$
onto $\mathbb{Z}_4$ by 
$\phi((1\: 3\: 5\: 7\: 9)(2\: 4\: 6\:8\:10))=0$,
and $\phi((1\: 2\: 9\: 8)(3\: 6\: 7\: 
4)(5\: 10))=1$.
\item $\Gamma'=
\mathrm{Ker}(\phi)=
\langle
(1\: 3\: 5\: 7\: 9)(2\: 4\: 6\: 8\: 10)
\rangle$.
\end{itemize}
Then
\begin{itemize}
\item $\Gamma' \triangleleft \Gamma$,
\item $\Gamma/\Gamma'=\mathbb{Z}_4$ is
cyclic.
\item $|\Gamma'|=5$ is a prime power.
\end{itemize}
The vertex set of $\Delta^\Gamma$ is 
$\{G_{\{1,3\}}, G_{\{2\}}, G_{\{4\}},
G_{\{5\}}\}$, so by theorem \ref{theorem:1},
\[\begin{array}{rcl}
\chi(\Delta^\Gamma) &=&
i_{13}+i_{2}+i_{4}+i_{5}
-i_{123}-i_{134}-i_{135}-i_{24}-i_{25}-i_{45}
+i_{1234}+i_{1235}+i_{1345}+i_{245}\\
&=&
i_{13}+2i_{2}+1
-2i_{123}-i_{135}-i_{24}-2i_{25}
+0+2i_{1235}+i_{245}\\
&=&2i_{2}+i_{13}-i_{24}-2i_{25}
-2i_{123}-i_{135}+i_{245}+2i_{1235}+1\\
&=&1\\
\end{array}\]
\end{proof}


\begin{lemma}\label{lemma:T6}
For any simplicial complex in $\mathfrak{C}$,
\[2i_2-i_{24}+i_{135}-2i_{1235}=1\]
\end{lemma}
\begin{proof}
Let
\begin{itemize}
\item $\Gamma=\langle
(2\: 4\: 6\: 8\: 10),
(1\: 6)(2\: 7)(3\: 8)(4\: 9)(5\: 10)\rangle$.
\item Let $\phi$ be a homomorphism from $\Gamma$
onto $\mathbb{Z}_{10}$ by
$\phi( (2\: 4\: 6\: 8\: 10) ) = 4$ and
$\phi( (1\: 6)(2\: 7)(3\: 8)(4\: 9)(5\: 10) )=5$.
\item $\Gamma'=
\mathrm{Ker}(\phi)=
\langle
(1\: 3\: 5\: 7\: 9)(2\: 10\: 8\: 6\: 4)
\rangle$.
\end{itemize}
Then
\begin{itemize}
\item $\Gamma' \triangleleft \Gamma$,
\item $\Gamma/\Gamma'=\mathbb{Z}_{10}$ is cyclic.
\item $|\Gamma'|=5$ is a prime power.
\end{itemize}
The vertex set of $\Delta^\Gamma$ is 
$\{G_{\{1,3,5\}}, G_{\{2\}}, G_{\{4\}}\}$, 
so by theorem \ref{theorem:1},
\[\begin{array}{rcl}
\chi(\Delta^\Gamma) &=&
i_{135}+i_{2}+i_{4}
-i_{24}-i_{1235}-i_{1345}\\
&=&
2i_{2}-i_{24}+i_{135}-2i_{1235}\\
&=&1\\
\end{array}\]
\end{proof}


\begin{lemma}\label{lemma:T8}
For any simplicial complex in $\mathfrak{C}$,
\[i_{14}=i_{145}\]
\end{lemma}
\begin{proof}
Let
\begin{itemize}
\item $\Gamma=\langle
(1\: 3\: 5\: 7\: 9)(2\: 4\: 6\: 8\: 10),
(2\: 7)(5\: 10)\rangle$.
\item Let $\phi$ be the homomorphism from 
$\Gamma$ onto $\mathbb{Z}_5$ defined by
$\phi((1\:3\:5\:7\:9)(2\:4\:6\:8\:10))=1$,
and $\phi((2\: 7)(5\: 10))=0$.
\item $\Gamma'=
\mathrm{Ker}(\phi)=
\langle
(1\: 6)(2\: 7),
(1\: 6)(5\: 10)
(3\: 8)(4\: 9),
(4\: 9)(5\: 10)
\rangle$.
\end{itemize}
Then
\begin{itemize}
\item $\Gamma' \triangleleft \Gamma$,
\item $\Gamma/\Gamma'=\mathbb{Z}_5$ is cyclic.
\item $|\Gamma'|=2^4$ is a prime power.
\end{itemize}
The vertex set of $\Delta^\Gamma$ is 
$\{G_{\{1,4\}}, G_{\{2,3\}}, G_{\{5\}}\}$, 
so by theorem \ref{theorem:1},
\[\begin{array}{rcl}
\chi(\Delta^\Gamma) &=&
i_{14}+i_{23}+i_{5}
-i_{145}-i_{235}-i_{1234}\\
&=& 2i_{14}+1-2i_{145}-0 \\
&=& 1\\
\end{array}\]
\end{proof}


This lemma can be directly deduced from lemma
\ref{lemma:T6}, but it is not clear which one
could be generalized most.
\begin{lemma}\label{lemma:T18}
For any simplicial complex in $\mathfrak{C}$,
\[i_{24}+i_{135}=1\]
\end{lemma}
\begin{proof}
Let
\begin{itemize}
\item $\Gamma=\langle
(1\: 7\: 9\: 3)(2\: 4\: 8\: 6),
(2\: 4\: 6\: 8\: 10),
(1\: 4\: 3\: 2\: 9\: 6\: 7\: 8)(5\: 10)
\rangle$.
\item Let $\phi$ be the homomorphism from 
$\Gamma$ onto $\mathbb{Z}_8$ by $\phi(
(1\: 7\: 9\: 3)(2\: 4\: 8\: 6))=2$,
$\phi((2\: 4\: 6\: 8\: 10))=0$, and
$\phi((1\: 4\: 3\: 2\: 9\: 6\: 7\: 8)
(5\: 10))=3$.
\item $\Gamma'= \mathrm{Ker}(\phi) = 
\langle
(2\: 4\: 6\: 8\: 10),
(1\: 5\: 9\: 3\: 7)(2\: 10\: 8\: 6\: 4)
\rangle$.
\end{itemize}
Then
\begin{itemize}
\item $\Gamma' \triangleleft \Gamma$,
\item $\Gamma/\Gamma'=\mathbb{Z}_8$ is cyclic.
\item $|\Gamma'|=5^2$ is a prime power.
\end{itemize}
The vertex set of $\Delta^\Gamma$ is 
$\{G_{\{1,3,5\}}, G_{\{2,4\}}\}$, 
so by theorem \ref{theorem:1},
\[\chi(\Delta^\Gamma) = i_{135}+i_{24} =1.\]
\end{proof}

\section{Gathering the facts}

The indicators $i_{3}, i_{4}, i_{23},
i_{34}, i_{35}, i_{45}, i_{134}, i_{234},
i_{235}, i_{345}, i_{1245}, i_{1234},$
and $i_{1235}$ are equal to another
indicator by graph isomorphism. By 
lemma \ref{lemma:T24}, $i_\emptyset=
i_5=1$ and $i_{1234}=i_{12345}=0$.
We have the equalities
\[\begin{array}{r|l}
\mathrm{Lemma} & \mathrm{Equality}\\
\hline
\mathrm{\ref{lemma:T1}} &
2i_{1}+2i_{2}-2i_{12}-i_{13}-
2i_{14}-2i_{15}-i_{24}-2i_{25}
+2i_{123}+\\
& 2i_{124}
+2i_{125}+i_{135}+2i_{145}
+ i_{245}-2i_{1235}-2i_{1245}=0\\
\hline
\mathrm{\ref{lemma:T4}} &
2i_2+i_{13}-i_{24}-2i_{25}-2i_{123}-i_{135}
+i_{245}+2i_{1235}=0\\
\hline
\mathrm{\ref{lemma:T6}} &
2i_2-i_{24}+i_{135}-2i_{1235}=1\\
\hline
\mathrm{\ref{lemma:T8}} &
i_{14}=i_{145} \\
\hline
\mathrm{\ref{lemma:T18}} &
i_{24}+i_{135}=1\\
\end{array}\]

\begin{theorem}\label{lemma:six}
For any simplicial complex in $\mathfrak{C}$,
the indicators of the transitive graphs are
as one of the following six columns.
\begin{center}
\emph{
\begin{tabular}{c|cccccc}
Indicator$\setminus$Name & 
$\mathbf{A}$ &
$\mathbf{A}^\ast$ &
$\mathbf{B}$ &
$\mathbf{B}^\ast$ &
$\mathbf{C}$ &
$\mathbf{C}^\ast$\\
\hline
$i_{1}$    & 1 & 1 & 1 & 1 & 1 & 1 \\
$i_{2}$    & 0 & 1 & 1 & 1 & 1 & 1 \\
$i_{5}$    & 1 & 1 & 1 & 1 & 1 & 1 \\
$i_{12}$   & 0 & 1 & 1 & 1 & 0 & 1 \\
$i_{13}$   & 1 & 1 & 0 & 0 & 0 & 0 \\
$i_{14}$   & 0 & 1 & 0 & 1 & 0 & 1 \\
$i_{15}$   & 1 & 1 & 0 & 1 & 1 & 1 \\
$i_{24}$   & 0 & 0 & 1 & 1 & 1 & 1 \\
$i_{25}$   & 0 & 1 & 1 & 1 & 1 & 1 \\
$i_{123}$  & 0 & 1 & 0 & 0 & 0 & 0 \\
$i_{124}$  & 0 & 0 & 0 & 1 & 0 & 0 \\
$i_{125}$  & 0 & 1 & 0 & 0 & 0 & 1 \\
$i_{135}$  & 1 & 1 & 0 & 0 & 0 & 0 \\
$i_{145}$  & 0 & 1 & 0 & 1 & 0 & 1 \\
$i_{245}$  & 0 & 0 & 1 & 1 & 1 & 1 \\
$i_{1234}$ & 0 & 0 & 0 & 0 & 0 & 0 \\
$i_{1235}$ & 0 & 1 & 0 & 0 & 0 & 0 \\
$i_{1245}$ & 0 & 0 & 0 & 0 & 0 & 0 \\
$i_{P}$    & 0 & 1 & 1 & 1 & 1 & 1 \\
$i_{\bar{P}}$ &0&1 & 0 & 0 & 0 & 0 \\
\end{tabular}}
\end{center}
\end{theorem}
\begin{proof}
The proof is in three cases.\\
\noindent \textbf{Case 1}: $i_2=0$.

The indicators of all graph including 
$G_{\{2\}}$ are zero. The only undetermined are
$i_1,i_{15},i_{13}$, and $i_{135}$. By lemma
\ref{lemma:T18}, $i_{24}+i_{135}=1\Rightarrow
i_{135}=1$. Since $i_1\geq i_{15} \geq i_{13}
\geq i_{135}$, they are all equal to 1.
Since $i_2=0$, $i_P=i_{\bar{P}}=0$. This
is column $\mathbf{A}$.\\
\noindent \textbf{Case 2}: $i_{1235}=1$.

This is the dual assumption of $i_2=0$, which
gives column $\mathbf{A}^\ast$.\\
\noindent \textbf{Case 3}: 
$i_2=1$ and $i_{1235}=0$.

Lemma \ref{lemma:T6} becomes $-i_{24}+i_{135}
=-1$, hence $i_{24}=1$ and $i_{135}=0$.
Lemma \ref{lemma:T4} becomes $i_{13}-2i_{25}
-2i_{123}+i_{245}=-1$. Both $i_{25}$ and
$i_{123}$ cannot be 1, and $i_{123}\leq i_{25}$,
so $i_{123}=0$. And now $i_{25}=1$ since 
$i_{13}-2i_{25}+i_{245}=-1$. Thus $i_{13}+
i_{245}=1$. Note that $i_{1}=1$ and $i_{1245}=0$
since $i_{25}=1$ and $i_{123}=0$.\\
\noindent \textbf{Case 3.1}: $i_{13}=1$ and
$i_{245}=0$.

From subgraph inclusion we get $i_{15}=1$ and
$i_{124}=0$. Remove $i_{14}$ and $i_{145}$ 
from lemma \ref{lemma:T1} by using lemma
\ref{lemma:T8}, and insert values to get
$-2i_{12}+2i_{124}=2$. But this can never be
true since $i_{12}\geq i_{124}$.\\
\textbf{Case 3.2}: $i_{13}=0$ and
$i_{245}=1$.

Once again removing $i_{14}$ and $i_{145}$ 
from lemma \ref{lemma:T1} by using lemma
\ref{lemma:T8}, and insert values, gives
$i_{12}+i_{15}-i_{124}-i_{125}=1$.
Since both $i_{12}$ and $i_{15}$ are greater
or equal with both $i_{124}$ and $i_{125}$, we 
have four different options, which are the
columns $\mathbf{B},\mathbf{B}^\ast,
\mathbf{C},\mathbf{C}^\ast$. The value of
$i_{14}=i_{145}$ is uniquely determined in
each column since $i_{12},i_{15}\geq i_{145}
=i_{14}\geq i_{124},i_{125}$. Finally
about the Petersen graph: $i_P\geq i_{245}=1$
and $i_{\bar{P}}\leq i_{13}=0$.
\end{proof}

\section{Computational aspects}
The proofs in this paper are computer
independent but to find the right lemmas for
section \ref{se} several algorithms were 
implemented in MAPLE. First a library of
subgroups of $S_{10}$ satisfying the conditions 
of Theorem~\ref{theorem:1} and \ref{theorem:2}
was constructed. This library, of what is called
non Oliver groups, was not complete. Each
group gives a linear equation of indicators, 
and only those with indicators of
transitive graphs were used. This system of
equations was heavily linearly dependent, so 
almost all equations with their groups were
thrown away. Those left became the lemmas of
section \ref{se}.

One interesting continuation of this
work is the removal of conditions on the 
indicators. That would give a huge equation
system, but maybe sufficient conditions on the
graph properties to do the remaining search
for a counterexample by brute force.

\end{document}